\documentclass[11pt]{article}
\usepackage{graphicx,mathrsfs,amssymb, color}
\usepackage{amsmath,latexsym,amsbsy,amssymb}
\usepackage{graphicx}

\def\ds{\displaystyle}
\def\forall{\hbox{for all}~}

\def\argmax{\hbox{arg}\!\max}

\def\ve{\varepsilon}
\def\n{\noindent}

\def\R{\mathbb{R}}

\def\vp{\varphi}

\def\vs{\vskip 2em}

\def\v{\vskip 1em}
\def\O{{\cal O}}

\def\bega{\begin{array}}
\def\enda{\end{array}}
\def\begi{\begin{itemize}}
\def\endi{\end{itemize}}

\def\Tilde{\widetilde}

\def\bel{\begin{equation}\label}
\def\eeq{\end{equation}}
\def\sqr#1#2{\vbox{\hrule height .#2pt
\hbox{\vrule width .#2pt height #1pt \kern #1pt
\vrule width .#2pt}\hrule height .#2pt }}
\def\square{\sqr74}
\def\endproof{\hphantom{MM}\hfill\llap{$\square$}\goodbreak}

\newtheorem{theorem}{Theorem}[section]
\newtheorem{corollary}[theorem]{Corollary}
\newtheorem{definition}[theorem]{Definition}
\newtheorem{remark}[theorem]{Remark}
\newtheorem{example}[theorem]{Example}
\newtheorem{lemma}[theorem]{Lemma}
\newtheorem{proposition}[theorem]{Proposition}

\begin{document}
\title{\bf On the structure  of the value function  \\of optimal exit time problems}\vs
\author{Piermarco Cannarsa$^{(1)}$, Marco Mazzola$^{(2)}$, and Khai T. Nguyen$^{(3)}$\\ 
\\
 {\small $^{(1)}$ Department of Mathematics, University of Roma ``Tor Vergata", }\\  
 {\small $^{(2)}$ Sorbonne Universit\'e and Universit\'e Paris Cit\'e, CNRS, IMJ-PRG, F-75005 Paris, France,}\\
  {\small $^{(3)}$ Department of Mathematics, North Carolina State University. }\\ 
 \\  {\small e-mails: cannarsa@axp.mat.uniroma2.it, ~marco.mazzola@imj-prg.fr,  ~ khai@math.ncsu.edu}
 \quad\\
 \quad\\
{\small (dedicated to Giovanni Colombo on the occasion of his 65th birthday)}
 }
\maketitle

\begin{abstract}  In this paper, we study  an optimal exit time problem   with  general running and terminal costs  and a target $\mathcal{S}\subset\R^d$  having an inner ball property  for a nonlinear control system that satisfies  mild  controllability assumptions. In particular,  Petrov's condition at the boundary of $\mathcal{S}$ is not required and  the value function $V$ may fail to be locally Lipschitz. In such a weakened set-up, we first establish a representation formula for proximal (horizontal) supergradients of $V$
by using transported proximal normal vectors. This allows us to obtain an external sphere condition for the hypograph of $V$ which yields several regularity properties. In particular, $V$ is almost everywhere twice differentiable and  the Hausdorff dimension of  its singularities is not greater than $d-1/2$. Furthermore,  besides optimality conditions for trajectories of the  optimal control problem,  we extend the analysis to propagation of singularities and  differentiability properties of the value function. An upper bound for the Hausdorff measure of the singular set  is also studied,  which implies that $V$  is a function of special bounded variation.
\\
		\quad\\
		{\bf Keyword.} Optimal exit time, proximal  supergradients, external sphere condition
		\medskip
		
		\n {\bf AMS Mathematics Subject Classification. 49N60, 49N05, 49J52, 49E30.} 
 \end{abstract}
 


\section{Introduction}
\label{sec:1} The present  paper is  concerned with  a rather general class of optimal exit time problems with a nonlinear control dynamics
 \bel{Cont}
 \begin{cases}
 y'(s)~=~f(y(s),u(s)),\\
 u(s)~\in~U,\\
 y(0)~=~x~\in~\R^d,
 \end{cases}
 \eeq
 where the function $f: \R^d\times U\to \R^d$ is  (for simplicity) $\mathcal{C}^2$, the control set $U$ is a nonempty compact subset of $\R^m$, and  $u(\cdot)$ belongs to   the set of admissible controls 
 \bel{Uad}
 \mathcal{U}_{ad}~\doteq~\left\{u:[0,\infty[\to U: u~\text{is~measurable}\right\}.
 \eeq
 Suppose we are now given a closed nonempty target set $\mathcal{S}\subseteq\R^d$, a running cost $r:\R^d\times U\to \R$ and a  terminal cost $g:\R^d\to\R$.  The value function at $x\in\R^d$ is defined by 
\bel{V1}
V(x)~\doteq~\inf_{u\in\mathcal{U}_{ad}}\left(\int_{0}^{\tau^{x,u}}r\left(y^{x,u}(s),u(s)\right)ds+g\left(y^{x,u}(\tau^{x,u})\right)\right)
\eeq
with $y^{x,u}(\cdot)$ being the unique Carath\'eodory solution of (\ref{Cont})  and $\tau^{x,u}$ being  the first  time at which $y^{x,u}(\cdot)$ reaches the target $\mathcal{S}$. By the dynamic programming principle, if the value function $V$ is continuous in $ \mathcal{R}\backslash\mathcal{S}$ then it is a viscosity solution of the associated Hamilton-Jacobi equation 
\bel{HJB}
\begin{cases}
\ds\max_{u\in U}\{-f(x,u)\cdot \nabla V(x)-r(x,u)\}~=~0,\qquad x\in \mathrm{int} (\mathcal{R}\backslash\mathcal{S}), \\[3mm]
V(x)~=~g(x), \qquad x\in\partial\mathcal{S},
\end{cases}
\eeq
with $\mathcal{R}$ being the controllable set of (\ref{Cont}) to the target $\mathcal{S}$. This fact can be used as an alternative characterization of $V$. To derive  optimality conditions on the trajectories of the system and, sometimes, an optimal feedback control, one requires some  suitable type of generalized gradient. The strongest regularity property one can expect for $V$---in fairly general cases---is semiconcavity. A semiconcavity result was obtained  in \cite{PCC} under suitable smoothness assumptions on $f,r,g, \mathcal{S},$ and  Petrov's condition, i.e., there exists $\mu>0$ such that for every $x\in \partial\mathcal{S}$ and $p$ proximal normal to $\mathcal{S}$ at $x$,
\bel{Petrov}
\min_{u\in U} p\cdot f(x,u)~\leq~-\mu |p|,
\eeq
which guarantees the local Lipschitz continuity of the associated minimum time function \cite{P}. Thus, together with  a compatibility condition on the terminal cost $g$ near to the boundary of the target $\mathcal{S}$,  Petrov's condition  yields the local Lipschitz continuity of  $V$. This result is then applied to establish  necessary optimality conditions, through the formulation of a suitable version of the Maximum Principle, and to study the singular set of $V$.

In several optimal exit time problems such as the minimum time problem, controllability assumptions weaker than Petrov's condition hold and   $V$ is continuous, but generally fails to be locally Lipschitz (see for example in \cite{M1,M2,M3}). Thus, a very natural question is to try to understand whether the structure of $V$ remains unchanged if in the above setting the controllability assumptions are weakened. The main purpose of this paper is to study the regularity of value function $V$  for the general system (\ref{Cont}), running cost $r$, and terminal cost $g$, assuming that $V$ is just continuous. In this case, the strongest  regularity of $V$  will be certainly weaker than semiconcavity. However, it is known that semiconcavity  amounts to  local Lipschitz continuity plus a uniform exterior sphere condition for the hypograph. Thus, one expects that  weakend assumptions may still guarantee some geometric properties of the hypograph (or epigraph) of $V$, which can in turn be used to derive  BV estimates for $V$ and $\nabla V$, the existence of Taylor's expansion of order two around almost every point, and bounds for the singular set of $V$ \cite{GA,N,ANV}. In the case of  time optimal control ($r \equiv 1, g\equiv 0$), for linear dynamics and a convex target the answer is positive:  it was proved in \cite{GAW} that the epigraph of the minimum time function satisfies a uniform exterior sphere condition or even has  positive reach \cite{F} (or $\varphi$-convexity \cite{GA}). Later on, this semiconvexity-type result was  analysed deeper in \cite{GN1,GN2}, but just for  linear  or linear affine control systems in dimension two. To remove this restriction, it is useful to look, first, at semiconcavity-type result for the time optimal control problem. In this case, several results beyond semiconcavity were obtained for  parameterized nonlinear control systems in \cite{GN,N} and  for differential inclusions in \cite{PK,PKA}, by deriving an external sphere condition for  the hypograph of the minimum time function under mild controllability assumptions. 
However, in this spirit, a study of  $V$ with general running and terminal costs is completely  lacking.

The first  purpose  of the present paper aims at addressing the above issue. In Theorem \ref{Main1}, we show that under some  standard assumptions on the dynamics and cost functionals, if $V$ is continuous and the target $\mathcal{S}$ satisfies an internal sphere condition, then the hypograph of $V$ satisfies an exterior sphere condition. The continuity of $V$ in $\mathcal{R}$ holds under a small time controllability  assumption and  a compatibility condition on the terminal cost $g$ (see  in Lemma \ref{Con-V}).  As a consequence, from   \cite[Theorem 1.3]{ANV}, one has that  outside a closed set of locally  finite $\mathcal{H}^{d-1}$ measure, the hypograph of $V$ has positive reach \cite{NV}. In particular, it is semiconcave and admits a second order Taylor expansion around almost every point. Moreover,  the Hausdorff dimension of  its singularities  is not greater than $d-1/2$, which is sharp (see in \cite[Proposition 7.3]{ANV}). In the meantime, we establish a new representation formula of {\it proximal supergradients of $V$} and {\it proximal horizontal supergradients of $V$} at every point $x\in \mathcal{R}\backslash\mathcal{S}$,
\bel{rela-1}
\partial^{P}V(x)~\supseteq~ \partial^{\infty}V(x)+\mathrm{co}\big[N^{\dagger}_1(x)\big],\qquad \partial^{\infty}V(x)~\supseteq~\mathrm{co}\big[ N^{\dagger}_0(x)\cup \{0\}\big],
\eeq
where $N^{\dagger}_0(x)$ and $N^{\dagger}_1(x)$ are the two sets of transported proximal normals introduced in  Definition \ref{TP}.
In general the inclusions in  (\ref{rela-1}) can be strict but we prove in Theorem \ref{Re-pre} that  equalities hold at every point $x\in \mathrm{int} (\mathcal{R}\backslash\mathcal{S})$ such that the proximal normal cone to the hypograph of $V$ is pointed, i.e.
\bel{pointed-1}
N^{P}_{\mathrm{hypo}(x)}(x,V(x))\bigcap\left(-N^{P}_{\mathrm{hypo}(x)}(x,V(x))\right)~=~\{0\}.
\eeq
The pointedness assumption plays a major role: actually, exposed rays of the normal cone to the hypograph are special normals, since  they can be approximated by normals at differentiability points of $V$, as established in Proposition \ref{Pres-f}.  In particular, if (\ref{pointed-1}) holds for every $x\in \mathrm{int} \mathcal{R}\backslash\mathcal{S}$, then the hypograph of $V$ has positive reach.



In the second part of the paper, as an application of the obtained results, we study  optimality conditions for the trajectories of our control problem. In Section \ref{sec4}, we give a complete version of  the Pontryagin's maximum principle for  both non-horizontal  and horizontal cases,  that generalizes similar results  to continuous value function and nonsmooth target for every proximal normal vector to the target.   Then, we  provide sufficient conditions for optimality and relate the uniqueness of optimal trajectories to  a differentiability property of the hypograph of $V$. In Section \ref{sec5}, under a suitable assumption on the boundary of the target $\mathcal{S}$, we show that the singular set of $V$ is countably $\mathcal{H}^{d-1}$-rectifiable. As a consequence, in this case $V$  is a function of special bounded variation in $ \mathcal{S}^c$. Finally, using the notion of optimal point, we establish  results on the propagation of singularities and  the differentiability of $V$ along optimal trajectories.


%
%
%
%
%

\section{Notation and preliminary results}
\setcounter{equation}{0}
\subsection{Nonsmooth analysis}
Given a closed set $Q\subseteq\R^d$, we say that a vector $v\in\R^d$ is a {\it Fr\'echet normal} to $Q$ at the point  $x\in \partial Q$ (and we write $v\in N^{F}_Q(x)$) if 
\[
\limsup_{Q\ni y\to x}~v\cdot {y-x\over |y-x|}~\leq~0.
\]
We say that is it a {\it proximal normal} to $Q$ at $x\in Q$ \big(and write $v\in N^{P}_{Q}(x)$\big) if there exists $\sigma\geq 0$ such that 
\bel{prox}
v\cdot (y-x)~\leq~\sigma\cdot |y-x|^2\qquad\forall y\in Q.
\eeq
Equivalently, $v\in N^{P}_{Q}(x)$ if and only if there is $\lambda>0$ such that $d_{Q}(x+\lambda v)=\lambda|v|$, with $d_Q(y)\doteq\ds\min_{x\in Q}\{|y-x|\}$ being the distance from $y$ to $Q$. We say that $v\in N^{P}_{Q}(x)$ is {\it realized by a ball of radius $\rho$} if (\ref{prox}) holds for $\sigma=|v|/(2\rho)$.

\begin{definition}[Exterior sphere condition]\label{sphere-cond} Let $\rho:\partial Q\to (0,+\infty)$ be continuous. We say that $Q$ satisfies the $\rho(\cdot)$-exterior sphere condition if and only if for every $x\in\partial Q$ there exists a unit vector $v_x\in N^{P}_{Q}(x)$ realized by a ball of radius $\rho(x)$, i.e.,
\[
 v_x\cdot (y-x)~\leq~{1\over 2\rho(x)}\cdot |y-x|^2\qquad\forall y\in Q.
\]
\end{definition}
A set $Q$ satisfies the $\rho(\cdot)$-interior  sphere condition if the closure of the complement of $Q$ satisfies  the $\rho(\cdot)$-exterior sphere condition. If $Q$
satisfies the $\rho(\cdot)$-interior sphere condition for some constant function $\rho(\cdot)\equiv\rho_0$, then we also say that $Q$ has the inner ball property (of radius $\rho_0$). On the other hand, sets  characterized by a stronger exterior sphere condition are called $\vp$-convex sets (proximally smooth sets \cite{C1}, prox-regular sets \cite{NS}, or sets with positive reach) and are defined as follows: 
\begin{definition}[$\vp$-convex sets] Let $\vp:\partial Q\to (0,+\infty)$ be continuous. We say that 
a closed set $Q$ is $\vp$-convex if, for any $x\in\partial{Q}$, every proximal unit  normal $v_x\in N^{P}_{Q}(x)$ is realized by a ball of radius $\vp(x)$.
\end{definition}
 $\vp$-convex sets enjoy, locally, several properties that are typical of convex sets. In particular, this is the case for the metric projection, which
is unique in a neighborhood of a $\vp$-convex set. This fact is used in proving all the
regularity properties which are satisfied by $\vp$-convex sets  (see, e.g., \cite[Section 4]{F}).  The following equivalence  result between the exterior sphere condition and
$\vp$-convexity was obtained in \cite{GN,NS} under a geometric condition.
\begin{proposition}\label{EQ} Let $Q\subseteq\R^d$ be a closed set satisfying  the $\rho(\cdot)$-exterior sphere condition. If $Q$ is pointed, i.e., 
\[
N^{P}_Q(x)\bigcap \left(-N^{P}_Q(x)\right)~=~\{0\}\qquad\forall x\in \partial{Q},
\]
then   $Q$ is a $\vp$-convex set. 
\end{proposition}

Using techniques of nonsmooth analysis and geometric measure theory, further properties were deduced in \cite{N, NV, ANV}. In particular, a set satisfying an exterior (or interior) sphere condition has locally finite perimeter. In addition, its boundary is locally a finite union of Lipschitz graphs \cite[Corollary 4.3]{ANV}.
\medskip

In what follows, let $\Omega$ be a subset of $\R^d$  and $f:\Omega\to\R$ be an upper semicontinuous function.  By using the hypograph of $f$, which  is denoted by 
\[
\mathrm{hypo}(f)~\doteq~\left\{(x,\beta)\in \Omega\times\R:\beta\leq f(x)\right\},
\]
we can define some generalized differentials for $f$.  More precisely, for every $x\in \Omega$ and $v\in \R^d$, we say that
\begin{itemize}
\item  $v$ is a {\it proximal supergradient} of $f$ at $x$ \big($v\in\partial^{P}f(x)$\big) if $(-v,1) \in N^P_{\mathrm{hypo}(f)}(x,f(x))$;  equivalently (see \cite[Theorem 1.2.5]{CLSW}), $v\in \partial^{P}f(x)$ if and only if there exist $\sigma,\eta>0$ such that 
\bel{Sub-D}
f(y)~\leq~f(x)+v\cdot (y-x)+\sigma\cdot |y-x|^2\qquad\forall y\in B(x,\eta)\cap \Omega;
\eeq
\item $v$ is a {\it proximal horizontal supergradient} of $f$ at $x$ \big($v\in\partial^{\infty}f(x)$\big) if $(-v,0) \in N^P_{\mathrm{hypo}(f)}(x,f(x))$, i.e., there exists $\sigma>0$ such that 
\bel{Sub-HD}
-v \cdot (y-x)~\leq~\sigma\cdot\left(|y-x|^2+|\beta-f(x)|^2\right)\qquad\forall (y,\beta)\in \mathrm{hypo}(f);
\eeq
\item $v$ is a {\it reachable gradient} of $f$ at $x$ \big($v\in\partial^{*}f(x)$\big) if there exists a sequence of differentiable points $x_n\in\Omega\backslash\{x\}$ of $f$, converging to $x$, such that  $\ds v=\lim_{n\to\infty}Df(x_n)$;
 \item $v\neq 0$ is a {\it horizontal reachable gradient} of $f$ at $x$ \big($v\in\partial^{*,\infty}f(x)$\big) if there exists a sequence of differentiable points $x_n\in\Omega\backslash\{x\}$ of $f$ converging to $x$ such that 
\[
{v\over |v|}~=~\lim_{n\to\infty} {Df(x_n)\over |Df(x_n)|}\qquad\mathrm{with}\qquad \lim_{n\to\infty}|Df(x_n)|=\infty.
\]
\end{itemize}
The proximal horizontal superdifferential $\partial^{\infty}f$  plays an important role in the study of a certain class of non-Lipschitz functions \cite{GA,N,ANV}.  More precisely, adapting Definition \ref{sphere-cond}, we introduce a class of upper semicontinuous functions  inheriting many regularity properties of semiconcave functions, which have several applications to both optimal control and partial differential equations  (see \cite{BD,PC}). 
\begin{definition}
We say that $\mathrm{hypo}(f)$  satisfies the $\rho(\cdot)$-exterior sphere condition if for every $x\in\Omega$ there exists a nonzero vector $v_x\in N^{P}_{\mathrm{hypo}(f)}(x,f(x))$ which is realized by a ball in  $\R^{d+1}$ of radius $\rho(x)$.
 \end{definition}
  It is well-known (see e.g. in \cite{ANV,N}) that every Lipschitz function $f:\Omega\to\R$ with $\mathrm{hypo}(f)$ satisfying a $\rho(\cdot)$-exterior sphere condition is {\it locally semiconcave}, i.e., for every $x\in \Omega$ there exist $\delta,c>0$ such that  for every $x_0,x_1\in B(x,\delta)$ and $\lambda\in [0,1]$, it holds 
 \[
 (1-\lambda)f(x_0)+\lambda f(x_1)-f(x_{\lambda})~\leq~c\cdot \lambda(1-\lambda)|x_0-x_1|^2,
 \]
where $x_\lambda\doteq(1-\lambda)x_0+\lambda x_1$. Equivalently, the map $z\mapsto f(z)-c\cdot |z|^2$ is concave in $B(x,\delta)$.
Thus, $f$ enjoys several properties of a concave function, in particular a.e. twice differentiability and the {\it superdifferential} of  $f$ at $x$ coincides with $\partial^{P}f(x)$, i.e., 
\[
\partial^{P}f(x)~=~D^+f(x)~\doteq~\left\{p\in \R^d:\underset{y\to x}{\lim\sup}\dfrac{f(y)-f(x)-p\cdot(y-x)}{|y-x|}~\leq~0\right\}.
\]
In the case where $f$ is just upper semicontinuous, let us denote by 
\bel{f-infty}
\Sigma_{f,\infty}~\doteq~\left\{x\in \Omega: \partial^{\infty}f(x)\neq \{0\}\right\}
\eeq
the set of points where the proximal horizontal superdifferential of $f$  contains a nonzero vector.  Here we recall a result which is a combination of the ones in   \cite{PKA,GA,ANV,N}.
\begin{proposition}\label{hypo-upp} Given an open set $\Omega\subset\R^d$, let $f:\Omega\to\R$ be upper semicontinuous such that  $\mathrm{hypo}(f)$ satisfies the $\rho(\cdot)$-exterior sphere condition. Then the function $f$ is in $BV_{loc}(\Omega)$ and   $\Sigma_{f,\infty}$ coincides with the set of all non-Lipschitz points of $f$. Moreover, $\Sigma_{f,\infty}$  is a closed subset  of $\Omega$ and the following holds true:
\begin{itemize}
\item [(i)]  the Hausdorff dimension of $\Sigma_{f,\infty}$ is not greater than $d-1/2$;
\item [(ii)] the function $f$ is locally semiconcave in the open set $\Omega\backslash\Sigma_{f,\infty}$.
\end{itemize}
Consequently, $f$ is a.e. Fr\'echet differentiable and admits a second order Taylor expansion around a.e. point in $\Omega\backslash\Sigma_{f,\infty}$. In addition, the set of points where the graph of f is nonsmooth has small Hausdorff dimension. More precisely, for every $k=1,2,\dots, d$, the set $\left\{x\in\Omega: the~dimension~of~\partial^{P}f(x)~is~\geq~k  \right\}$ is countably $\mathcal{H}^{d-k}$-rectifiable. 
\end{proposition}

In addition, let us establish necessary and sufficient conditions for the differentiability of the hypograph of $f$ at a point $(x,f(x))$ which is defined as follows: 
\begin{definition} Let $x\in\Omega$. We say that  $\mathrm{hypo}(f)$ is differentiable at $(x,f(x))$ if
$$N^{F}_{{\mathrm{hypo}(f)}}(x,f(x))~=~\{s v_x:s\in [0,\infty[\},$$ 
for some unit vector $v_x\in \R^{d+1}$ satisfying 
\bel{con-dd1}
v_x\cdot (y-x,f(y)-f(x))~=~o\big(|y-x|+|f(y)-f(x)|\big)\qquad\forall y\in\Omega.
\eeq
\end{definition}

\begin{proposition}\label{Diff} 
Under the same assumptions of Proposition \ref{hypo-upp}, $\mathrm{hypo}(f)$ is differentiable at $(x,f(x))$ if and only if 
\bel{dim-1}
N^{P}_{\mathrm{hypo}(f)}(x,f(x))~=~\{s v_x: s\in [0,\infty[\}
\eeq
for some unit vector $v_x\in \R^{d+1}$.
\end{proposition}
{\bf Proof.} {\bf 1.} Assume that $\mathrm{hypo}(f)$ is differentiable at $(x,f(x))$. Then by the definition we have 
\[
N^{P}_{\mathrm{hypo}(f)}(x,f(x))~\subseteq~N^{F}_{\mathrm{hypo}(f)}(x,f(x))~=~\{s v_x:s\in [0,\infty[\},
\]
for some  unit vector $v_x\in \R^{d+1}$.  By the exterior sphere condition on $\mathrm{hypo}(f)$, we have that $N^{P}_{\mathrm{hypo}(f)}(x,f(x))\backslash \{0\}$ is non-empty, and this yields  (\ref{dim-1}).
\medskip

{\bf 2.}  On the other hand if  (\ref{dim-1}) holds for  some unit vector $v_x\in \R^{d+1}$ then 
\bel{dqc}
v_x\cdot (y-x,\beta-f(x))~\leq~\O(1)\cdot \left(|y-x|^2+|\beta-f(x)|^2\right)\quad\forall (y,\beta)\in \mathrm{hypo}(f).
\eeq
For every $y\in \Omega$, there exists a unit vector $v_y\in N^{P}_{\mathrm{hypo}(f)}(y,f(y))$ realized by a ball of radius $\rho(y)$. In particular, one has 
 \[
v_y\cdot (x-y,f(x)-f(y))~\leq~{1\over 2\rho(y)}\cdot \left(|x-y|^2+|f(x)-f(y)|^2\right),
\]
and this implies 
\begin{multline*}
v_x\cdot \big(x-y,f(x)-f(y)\big)\\~\leq~\O(1)\cdot \big(|v_x-v_y|+|y-x|+|f(y)-f(x)|\big)\cdot \big(|y-x|+|f(y)-f(x)|\big).
\end{multline*}
Since $v_x$ is the unique unit vector $v_y\in N^{P}_{\mathrm{hypo}(f)}(x,f(x))$, one gets from   the $\rho(\cdot)$-exterior sphere condition of ${\mathrm{hypo}(f)}$ that $\ds\lim_{y\to x} |v_y-v_x|=0$. Thus, the above estimate and (\ref{dqc}) yield (\ref{con-dd1}).

{\bf 3.} Finally, to obtain the differentiability of ${\mathrm{hypo}(f)}$ at $(x,f(x))$, we need to show that 
\[
N^{F}_{{\mathrm{hypo}(f)}}(x,f(x))~=~\{s v_x:s\in [0,\infty[\}.
\]
Two cases are considered:
\medskip

 $\bullet$ If $v_x=\ds{(-\xi_x,1)\over |(-\xi_x,1)|}$ for some $\xi_x\in \R^d$, then $x\notin \Sigma_{f,\infty}$ and, by Proposition \ref{hypo-upp}, $f$ is semiconcave in a neighborhood of $x$. Thus, $f$ is differentiable at $x$ and 
\[
N^{F}_{{\mathrm{hypo}(f)}}(x,f(x))=\{s (-Df(x),1):s\in [0,\infty[\}.
\]
$\bullet$ Otherwise,  if $v_x=\ds {(-\xi_x,0)\over |(-\xi_x,0)|}$ for some $0\neq \xi_x\in \R^d$ then 
\bel{s-xix}
-{\xi_x\over |\xi_x|}\cdot (y-x)~\leq~{1\over 2\rho(x)}\cdot\left(|y-x|+|\beta-f(x)|\right)\quad\forall (y,\beta)\in \mathrm{hypo}(f).
\eeq
Suppose by contradiction that there exists a nonzero vector $(-\xi_1,\lambda_1)\in N^{F}_{\mathrm{hypo}(f)}(x,f(x))$ with $\lambda_1\in \{0,1\}$ and  $\ds {(-\xi_1,\lambda_1)\over |(-\xi_1,\lambda_1)|}\neq v_x$. Two subcases may occur:
\begin{itemize}
\item [-] If $\xi_1\notin[0,\infty[\,\cdot\,\xi_x$,
 then there is a unit vector $\zeta\in\R^d$ such that 
\bel{d-con1}
\zeta\cdot \xi_x~>~0,\qquad \zeta\cdot \xi_1~<~0.
\eeq
 By Proposition \ref{hypo-upp},  there exists a sequence $y_n\in\R^d$  such that $f$ is differentiable at $y_n$ and $ |y_n-(x+\zeta/n)|\leq n^{-2}$. Since $(-\xi_1,\lambda_1)\in N^{F}_{\mathrm{hypo}(f)}(x,f(x))$, one has 
 \[
 \begin{split}
 \lambda_1\cdot \big(\beta-f(x)\big)&~\leq~\xi_1\cdot (y_n-x)+o\big(|y_n-x|+|\beta-f(x)|\big)\\
 &~\leq~ {\xi_1\cdot\zeta\over n }+o\big(1/n+|\beta-f(x)|\big)\qquad\forall \beta\leq f(y_n),
 \end{split}
 \]
 and the second inequality in (\ref{d-con1}) yields 
 \[
 f(y_n)~<~f(x)\qquad\forall n>0~\mathrm{sufficiently~large}.
 \]
Moreover, by the exterior sphere condition of $\mathrm{hypo}(f)$, the vector $(-Df(y_n),1)\in  N^{P}_{\mathrm{hypo}(f)}(x,f(x))$ satisfies
\bel{ac2}
{(-Df(y_n),1)\over |(-Df(y_n),1)|}\cdot (z-y_n,\beta-f(y_n))~\leq~{1\over 2\rho(y_n)}\cdot \left(|z-y_n|^2+|\beta-f(y_n)|^2\right)
\eeq
for all $(z,\beta)\in \mathrm{hypo}(f)$. From  (\ref{dim-1}), we can suppose without loss of generality that $\lim_{n\to\infty} {(-Df(y_n),1)\over |(-Df(y_n),1)|}=~v_x$. For $n\geq 1$ sufficiently large, we can  choose $z=x$ and $\beta=f(y_n)$ in (\ref{ac2}) to derive 
\[
{(-Df(y_n),1)\over |(-Df(y_n),1)|}\cdot \left((x-y_n) n,0\right)~\leq~{n\over 2\rho(y_n)}\cdot |x-y_n|^2.
\]
Taking $n\to+\infty$, we obtain 
\bel{dda}
\xi_x\cdot \zeta~\leq~0
\eeq
and this contradicts the first inequality in (\ref{d-con1}).
\medskip

\item [-] Otherwise, if $\xi_1=s_1\xi_x$ for some $s_1\in [0,\infty[$, then $\lambda_1=1$ and 
\bel{s1s}
-s_1\xi_x\cdot (y-x)+\beta-f(x)~\leq~o(|y-x|+|\beta-f(x)|)
\eeq
for all $(y,\beta)\in \mathrm{hypo}(f)$. In particular, for some $\delta>0$ sufficiently small and $C>0$, 
\[
f(x)-f(y)~\geq~ -C\cdot |y-x|\qquad\forall y\in B(x,\delta).
\]
%
%
Let $\zeta\in\R^d$ be such that $\zeta\cdot \xi_x>0$ and let $(y_n)_n$ be a sequence in $B(x,\delta)$ such that $|y_n-(x+\zeta/n)|\leq n^{-2}$, $f$ is differentiable at $y_n$, (\ref{ac2}) holds for all $(z,\beta)\in \mathrm{hypo}(f)$, and 
$$
 \lim_{n\to\infty} {(-Df(y_n),1)\over |(-Df(y_n),1)|}~=~v_x=\ds{(-\xi_x,0)\over |(-\xi_x,0)|}.
$$
Choosing $z=x$ and $\beta=f(y_n)-C\cdot |y_n-x|$ in (\ref{ac2}) and multiplying both sides by $n$, we have
\[
n\cdot {(-Df(y_n),1)\over |(-Df(y_n),1)|}\cdot (x-y_n, -C\cdot |y_n-x|)~\leq~n\cdot {C^2+1\over 2\rho(y_n)}\cdot |x-y_n|^2.
\]
Taking $n\to+\infty$, we again obtain (\ref{dda}) and this yields a contradiction.
%
%
\end{itemize}

%
The proof is complete.
\endproof
\medskip
\medskip

To conclude this section, we recall some basic concepts concerning cones $C\subseteq \R^{d+1}$ \big(i.e., sets such that $\lambda v\in C$ for all $v\in C,\lambda \geq 0$\big). Then give a representation formula of $N^{P}_{\mathrm{hypo}(V)}(x,f(x))$ provided that $\mathrm{hypo}(f)$ satisfies an exterior sphere condition and $N^{P}_{\mathrm{hypo}(V)}(x,f(x))$ is pointed. Here  we say that C is {\it pointed} if 
$$C\cap \{-C\}~=~\{0\}.$$
\begin{definition}\label{expose}
A half line $\R^+\zeta\subseteq C$  is an exposed ray of a convex cone $C$ if there exists a unit vector $\bar{v}\in \R^d$ such that 
\[
\zeta\cdot \bar{v}~=~0\qquad\mathrm{and}\qquad w\cdot \bar{v}~<~0\quad\forall w\in C\backslash\big(\R^+\zeta\cup\{0\}\big).
\]
\end{definition}
It is known from \cite[Corollary 18.7.1, p.169]{R} that under the pointedness condition, a convex and closed cone $C$ can be represented by  its exposed rays. 
\begin{lemma}\label{co-exp} If a cone $\{0\}\neq C\subseteq \R^{d+1}$ is  closed, convex and pointed, then it is the closed convex hull of its exposed rays.
\end{lemma}

By Definition \ref{expose} and Lemma \ref{co-exp}, we show that  every exposed ray of $N^{P}_{\mathrm{hypo}(f)}(x,f(x))$ is generated by a reachable gradient or a horizontal reachable gradient of $f$ at $x$. Moreover, proximal (resp. proximal horizontal) supergradients can be computed in terms of reachable (resp. horizontal reachable) gradients.

\begin{proposition}\label{Pres-f} Given an open set $\Omega\subseteq\R^d$, let $f:\Omega\to \R$ be  continuous. Assume that $\mathrm{hypo}(f)$ satisfies the $\rho(\cdot)$-exterior sphere condition for some continuous function $\rho(\cdot)>0$. For every $x\in \Omega$, the followings hold:
\begin{itemize}
\item [(i).] The set of reachable gradients $\partial^{*}f(x)\bigcup \partial^{*,\infty}f(x)$ is non-empty and 
\bel{IN-CR}(-\partial^*f(x),1)\bigcup (-\partial^{*,\infty}f(x),0)\subseteq N^{P}_{\mathrm{hypo}(f)}(x,f(x)).
\eeq
\item [(ii).] Given any convex closed cone $C(x)$ such that  
$$(-\partial^*f(x),1)\bigcup (-\partial^{*,\infty}f(x),0)\subseteq C(x)\subseteq N^{P}_{\mathrm{hypo}(f)}(x,f(x)),$$ if a unit vector $(-\xi,\lambda)\in \R^{d}\times [0,\infty)$ belongs to an exposed ray of $C(x)$ then 
\bel{re1}
(-\xi,\lambda)~\in~\left\{{\zeta\over |\zeta|}:\zeta\in (-\partial^*f(x),1)\right\}\bigcup (-\partial^{*,\infty}f(x),0).
\eeq
\item [(iii).] If $N^{P}_{\mathrm{hypo}(f)}(x,f(x))$ is pointed then it is closed and 
\bel{hypo-p}
N^{P}_{\mathrm{hypo}(f)}(x,f(x))~=~\mathrm{co}\left\{\lambda v: \lambda\in [0,\infty), v\in  (-\partial^{*,\infty}f(x),0)\bigcup (-\partial^{*}f(x),1)\right\}.
\eeq
Moreover, the proximal (horizontal) supergradients of $f$ at $x$ are computed by 
\bel{sup-com}
\partial^{\infty}f(x)~=~\mathrm{co}\big[\partial^{*,\infty}f(x)\cup \{0\}\big],\qquad \partial^{P}f(x)~=~\mathrm{co}\big[\partial^{*}f(x)\big]+\partial^{\infty}f(x).
\eeq
\end{itemize}
\end{proposition}
{\bf Proof.} {\bf 1.} By Proposition \ref{hypo-upp}, the function $f$ is differentiable almost everywhere and this implies that $\partial^{*}f(x)\bigcup \partial^{*,\infty}f(x)$ is non empty. For every  sequence of differentiable point $x_n$ of $f$ converging to $x$, one has  that $(-Df(x_n),1)\in  N^{P}_{\mathrm{hypo}(f)}(x_n,f(x_n))$ is realized by a ball of radius $\rho(x_n)$, i.e., for all $(z,\beta)\in \mathrm{hypo}(f)$ it holds
\bel{Df-n}
{(-Df(x_n),1)\over |(-Df(x_n),1)|}\cdot (z-x_n,\beta-f(x_n))~\leq~{1\over 2\rho(x_n)}\cdot \left(|z-x_n|^2+|\beta-f(x_n)|^2\right).
\eeq
$\bullet$ If $\lim_{n\to\infty} Df(x_n)=\xi$ then taking $n\to\infty$ in (\ref{Df-n}), we get
\[
{(-\xi,1)\over |(-\xi,1)|}\cdot (z-x,\beta-f(x))~\leq~{1\over 2\rho(x)}\cdot \left(|z-x|^2+|\beta-f(x)|^2\right),
\]
and this yields $(-\xi,1)\in N^{P}_{\mathrm{hypo}(f)}(x,f(x))$.
\medskip

$\bullet$ If $\ds\lim_{n\to\infty} Df(x_n)=+\infty$ and $\ds\lim_{n\to\infty} {Df(x_n)\over |Df(x_n)|}=\xi$ then taking $n\to\infty$ in (\ref{Df-n}), we obtain
\[
-\xi\cdot (z-x)~\leq~{1\over 2\rho(x)}\cdot \left(|z-x|^2+|\beta-f(x)|^2\right)
\]
and this yields $(-\xi,0)\in N^{P}_{\mathrm{hypo}(f)}(x,f(x))$. Thus, the inclusion (\ref{IN-CR}) holds.
\medskip

{\bf 2.} Let  $(-\xi,\lambda)\in \R^{d}\times [0,\infty)$ be a unit vector which belongs to an exposed ray of $C(x)$. By Definition \ref{expose}, there exists $(v_0,\gamma)\in \R^d\times \R$ with unit norm such that for all $w\in C(x)\backslash \big(\R^+(-\xi,\lambda)\cup\{0\}\big)$
\bel{ex-sap}
(-\xi,\lambda)\cdot (v_0,\gamma)~=~0,\qquad w\cdot (v_0,\gamma)~<~0.
\eeq
By Proposition \ref{hypo-upp}, there exists a sequence $(x_n)_{n\geq 1}$ such that $f$ is differentiable at $x_n$ and 
\bel{x-n-c}
\begin{cases}
 \ds\lim_{n\to\infty}{(-Df(x_n),1)\over |(-Df(x_n),1)|}~=~(-\bar{\xi},\bar{\lambda})~\in~C(x),\\[4mm]
 \ds  \left|x+{v_0\over n}-{\xi\over n^{3/2}}-x_n\right|~\leq~{1\over n^2}\qquad\forall n\geq 1,
 \end{cases}
\eeq
From (\ref{ex-sap}) and (\ref{x-n-c}), one has $(-\bar{\xi},\bar{\lambda})\cdot (v_0,\gamma)\leq0$. Thus, in order to show that $(-\bar{\xi},\bar{\lambda})=(-\xi,\lambda)$ and obtain (\ref{re1}),
%
we only need to prove that 
\bel{k1}
(-\bar{\xi},\bar{\lambda})\cdot (v_0,\gamma)~\geq~0.
\eeq
Two cases are considered:
\medskip

{\bf Case 1.} Assume that $\lambda =0$. In this case, we have that  $\xi\cdot v_0=0$, $|\xi|=1$ and $(-\xi,0)\in C(x)\subseteq N^{P}_{\mathrm{hypo}(f)}(x,f(x))$ implies that for all $(z,\beta)\in \mathrm{hypo}(f)$
\[
 -\xi\cdot (z-x)~\leq~\O(1)\cdot \left(|z-x|^2+|\beta-f(x)|^2\right).
\]
 In particular, taking $z=x_n$ with $n$ sufficiently large, from the second inequality of (\ref{x-n-c})  we derive 
\[
{1\over n^{3/2}}-{1\over n^2}~\leq~\O(1)\cdot \left({1\over n^2}+|\beta-f(x)|^2\right)\quad\forall \beta\leq f(x_n),
\]
and this yields  
\[
f(x)-f(x_n)~=~|f(x)-f(x_n)|~\geq~{\O(1)\over n^{3/4}}.
\]
Hence, taking $z=x$ and $\beta=f(x)$ in (\ref{Df-n}), we have
\[
{Df(x_n)\over  |(-Df(x_n),1)|}\cdot {v_0\over n}~\leq~ \left({f(x)-f(x_n)\over 2\rho(x_n)}-{1\over |(-Df(x_n),1)|}\right)\cdot \big(f(x)-f(x_n)\big)+{\O(1)\over n^{3/2}}.
\]
This implies that $(-Df(x_n))_{n\geq 1}$ is unbounded, so that  $\bar{\lambda}=0$ and
\[
\begin{split}
(-\bar{\xi},\bar{\lambda})\cdot &(v_0,\gamma)~=~-\bar{\xi}\cdot v_0~=~\lim_{n\to\infty}{nDf(x_n)\over  |(-Df(x_n),1)|}\cdot (x-x_n)~\geq~-\lim_{n\to\infty}{n|x-x_n|^2\over 2\rho(x_n)}~=~0.
\end{split}
\]

{\bf Case 2.} Suppose that $\lambda >0$. In the case, we can assume $|v_0|=1$. Since  $(-{\xi\over \lambda},1)\in C(x)$, then
\bel{e-35}
-{\xi\over \lambda}\cdot (x_n-x)+ f(x_n)-f(x)~\leq~\O(1)\cdot \left(|x_n-x|^2+|f(x_n)-f(x)|^2\right)
\eeq
and this implies 
\bel{e-34}
f(x_n)-f(x)~\leq~\O(1)\cdot |x_n-x|\qquad\forall n\text{ large enough}.
\eeq
Thus, choosing $z=x$ and $\beta=f(x_n)-\O(1)\cdot |x_n-x|\leq f(x)$ in (\ref{Df-n}), we get
\[
{-Df(x_n)\over |(-Df(x_n),1)|}\cdot (x-x_n)~\leq~{\O(1)\cdot |x-x_n|\over  |(-Df(x_n),1)|}+{\O(1)\over 2\rho(x_n)}\cdot |x-x_n|^2,
\]
and (\ref{x-n-c}) yields 
\[
{-Df(x_n)\over |(-Df(x_n),1)|}\cdot v_0~\geq~-\O(1)\cdot\left({1\over |(-Df(x_n),1)|}+{1\over n^{1/2}}\right).
\]
If the sequence $(Df(x_n))_{n\geq 1}$ is unbounded then (\ref{ex-sap}) yields $\bar{\lambda}=0$. Taking $n\to\infty$ in the above estimate, we derive (\ref{k1}) by
\[
(-\bar{\xi},\bar{\lambda})\cdot (v_0,\gamma)~=~-\bar{\xi}\cdot v_0~=~\lim_{n\to\infty}{-Df(x_n)\over |(-Df(x_n),1)|}\cdot v_0~\geq~0.
\]
Otherwise, if  the sequence $(Df(x_n))_{n\geq 1}$ is bounded then  (\ref{x-n-c}) yields $\ds\lim_{n\to\infty}Df(x_n)={\bar{\xi}\over \bar{\lambda}}$ and $\bar{\lambda}>0$. Moreover, taking $z=x$ and $\beta=f(x)$ in (\ref{Df-n}), we have 
 \bel{e-36}
 -Df(x_n)\cdot (x-x_n)+f(x)-f(x_n)~\leq~{\O(1)\over 2\rho(x_n)}\cdot \left(|x-x_n|^2+|f(x)-f(x_n)|^2\right)
 \eeq
and (\ref{e-34}) implies  
\[
|f(x_n)-f(x)|~\leq~\O(1)\cdot |x_n-x|\qquad\forall n\geq 1.
\]
Summing (\ref{e-35}) and (\ref{e-36}), we have
\[
\left(-{\xi\over \lambda}+Df(x_n)\right)\cdot {x_n-x\over |x_n-x|}~\leq~\O(1)\cdot |x_n-x|.
\]
Taking $n\to\infty$, we obtain
\[
\left(-{\xi\over \lambda}+{\bar{\xi}\over \bar{\lambda}}\right)\cdot v_0~\leq~0\qquad\Longrightarrow\qquad \left(-{\xi\over \lambda},1\right)\cdot (v_0,\gamma)-\left(-{\bar{\xi}\over \bar{\lambda} },1\right)\cdot (v_0,\gamma)~\leq~0,
\]
and the first inequality of (\ref{ex-sap}) yields (\ref{k1}).
\medskip

{\bf 3.} Assume that $N^{P}_{\mathrm{hypo}(f)}(x,f(x))$ is pointed. Consider the convex cone 
\[
N_f(x)~\doteq~\mathrm{co}\left\{\lambda v: \lambda\in [0,\infty), v\in  (-\partial^{*,\infty}f(x),0)\bigcup (-\partial^{*}f(x),1)\right\}.
\]
Since $\mathrm{hypo}(f)$ satisfies the $\rho(\cdot)$-exterior sphere condition, one can easily show that  both $(-\partial^{*,\infty}f(x),0)$ and $(-\partial^{*}f(x),1)$ are closed subsets of $N^{P}_{\mathrm{hypo}(f)}(x,f(x))$. Thus, $N_f(x)$ is a convex, closed and pointed cone such that
\bel{sub-Nf}
N_f(x)~\subseteq~N^{P}_{\mathrm{hypo}(f)}(x,f(x)).
\eeq
We claim that the equality in (\ref{sub-Nf}) holds. Suppose by contradiction that there exists $\zeta\in N^{P}_{\mathrm{hypo}(f)}(x,f(x))\backslash N_f(x)$. The following cone 
\[
\Tilde{N}_f(x)~\doteq~\mathrm{co}\big[N_f(x)\cup \{\lambda \zeta: \lambda\in [0,\infty)\}\big]~\subseteq~N^{P}_{\mathrm{hypo}(f)}(x,f(x))
\]
 is  closed, convex, and pointed.  From (ii), all exposed rays of $\Tilde{N}_f(x)$ are in $N_f(x)$. Thus, by Lemma \ref{co-exp} we have 
 \[
 \Tilde{N}_f(x)~\subseteq~N_f(x)\,,
 \]
 and this yields a contradiction. 
 \medskip
 
 {\bf 4.} Finally, from (\ref{hypo-p}), one first has  
 \bel{iic}
 (-\partial^{*,\infty}f(x),0)~\subseteq~ \big(-\partial^{\infty}f(x),0\big)~\subseteq~\mathrm{co}\left\{\lambda v: \lambda\geq 0, v\in  (-\partial^{*,\infty}f(x),0)\right\},
 \eeq
  \bel{iic1}
\mathrm{co}\big[(-\partial^{*}f(x),1)\big]~\subseteq~ \big(-\partial^{P}f(x),1\big)~\subseteq~\mathrm{co}\left\{\lambda v: \lambda\geq 0, v\in  (-\partial^{*,\infty}f(x),0)\bigcup (-\partial^{*}f(x),1)\right\},
 \eeq
 and (\ref{iic}) implies that 
  \bel{conc1}
  \partial^{\infty}f(x)~=~\mathrm{co}\big\{\lambda v: \lambda\in [0,\infty), v\in \partial^{*,\infty}f(x) \big\}.
 \eeq
From the second inclusion in (\ref{iic1}), one has 
\[
\big(-\partial^{P}f(x),1\big)~\subseteq~\mathrm{co}\left\{\lambda v: \lambda\geq 0, v\in  (-\partial^{*,\infty}f(x),0)\right\}+\mathrm{co}\big[(-\partial^{*}f(x),1)\big],
\] 
and (\ref{conc1}) and the first inclusion of (\ref{iic1}) yields
\[
\partial^{\infty}f(x)+\mathrm{co} [\partial f^*(x)]~\subseteq~\partial^{P}f(x)~\subseteq~ \partial^{\infty}f(x)+\mathrm{co} [\partial f^*(x)].
\]
The proof is complete.
\endproof

\begin{corollary}\label{1-d-normal}
Under the same setting in Proposition \ref{Pres-f}, if $\mathrm{hypo}(f)$ is non differentiable at $(x,f(x))$ then the set of reachable gradients  $\partial^{*}f(x)\bigcup \partial^{*,\infty}f(x)$ contains  at least two elements. If in addition $\partial^{*}f(x)=\varnothing$, then $\partial^{*,\infty}f(x)$ contains two elements that are not positively collinear.
\end{corollary}
{\bf Proof.} Assume that $\mathrm{hypo}(f)$ is non differentiable at $(x,f(x))$. If $N^{P}_{\mathrm{hypo}(f)}(x,f(x))$ is pointed then from  (\ref{hypo-p})  and Proposition \ref{Diff}, the set of reachable gradients  $\partial^{*}f(x)\bigcup \partial^{*,\infty}f(x)$ contains  at least two elements  satisfying the required properties. Otherwise, if $N^{P}_{\mathrm{hypo}(f)}(x,f(x))$ is not pointed then $\{v_1, -v_1\}\subset\partial^{\infty}f(x)$ for some unit vector $v_1\in\R^d$ such that 
\bel{v1-1}
|v_1\cdot (z-x)|~\leq~\O(1)\cdot \left(|z-x|^2+|\beta-f(x)|^2\right)\quad\forall (z,\beta)\in\mathrm{hypo}(f).
\eeq
By (i) in Proposition \ref{Pres-f} and Proposition \ref{hypo-upp}, we can pick a reachable gradient $v^*\in \partial^{*}f(x)\bigcup \partial^{*,\infty}f(x)$ and  a sequence $(x_n)_{n\geq 1}$  such that $f$ is differentiable at $x_n$ and 
\[
\left|x+{v^*\over n}-{v_1\over n^{3/2}}-x_n\right|~\leq~{1\over n^2}\qquad\forall n\geq 1.
\]
Taking $z=x_n$ in \eqref{v1-1}, we obtain
\[
\left|{v_1v^*\over n}-{1\over n^{3/2}}\right|~\leq~\O(1)\cdot \left({1\over n^2}+|\beta-f(x)|^2\right)\qquad\forall\beta\leq f(x_n),
\]
and this yields 
\[
f(x)-f(x_n)~=~|f(x)-f(x_n)|~\geq~{O(1)\over n^{3/4}}\qquad \forall n~\text{large enough}.
\]
As in the proof of Proposition \ref{Pres-f}, $(-Df(x_n))_{n\geq 1}$ is unbounded so that  $\bar{\xi}\doteq\ds \lim_{n\to\infty}{Df(x_n)\over |Df(x_n)|}\in \partial^{*,\infty}f(x)$ and for all $(z,\beta)\in\mathrm{hypo}(f)$ it holds
\[
{(-Df(x_n),1)\over |(-Df(x_n),1)|}\cdot (z-x_n,\beta-f(x_n))~\leq~{1\over 2\rho(x_n)}\cdot \left(|z-x_n|^2+|\beta-f(x_n)|^2\right).
\]
Choosing $z=x$ and $\beta=f(x_n)$, we derive 
\[
\begin{split}
\bar{\xi}\cdot v^*&~=~\lim_{n\to\infty} \left[{- nDf(x_n)\over |(-Df(x_n),1)|}\cdot (x-x_n)\right]~\leq~\lim_{n\to\infty} {n |x-x_n|^2\over 2\rho(x_n) }~=~0.
\end{split}
\]
Thus, $\bar{\xi}\neq v^*$ and this completes the proof.
\endproof

\subsection{Optimal control}
In this subsection, we recall some basic concepts and well-known results on optimal exit time problems which will be used later. As usual, by the Landau symbol $\O(1)$ we shall denote a uniformly bounded quantity. Throughout the paper, we shall consider  the following standard hypothesis:
\begin{itemize}
\item [{\bf (H0).}] For every $x\in \R^d$, the set 
\[
\mathcal{F}(x)~\doteq~\left\{(v,\lambda)\in\R^d\times\R: \exists \omega\in U~~\text{such that}~~ v=f(x,\omega),\, \lambda \geq r(x,\omega)\right\}
\]
is convex.
\item [{\bf (H1).}] The map $f:\R^d\times U\to\R^d$ is continuous, uniformly bounded by a given constant $N$ such that $|f(x,u)|<N$, and satisfies  a  Lipschitz condition
\[
\sup_{\omega\in U}|f(y,\omega)-f(x,\omega)|~\leq~ \O(1)\cdot |y-x|\quad  \forall x,y\in\R^d.
\]
Moreover, the differential of $f$ with respect to the $x$ variable, denoted by $D_xf$, exists everywhere, is continuous with respect to both $x$ and $u$ and satisfies 
\[
\sup_{\omega\in U}|D_xf(y,\omega)-D_xf(x,\omega)|~\leq~ \O(1)\cdot |y-x|\quad\forall x,y\in\R^d.
\]
\item [{\bf (H2).}] The function $r:\R^d\times U\to\R$ is continuous and differentiable with respect to the $x$ variable. Both $r$ and  $D_xr$ are continuous and uniformly Lipschitz w.r.t $x$ such that 
\[
\sup_{\omega\in U}\big(|r(y,\omega)-r(x,\omega)|+\big|D_xr(y,\omega)-D_xr(x,\omega)\big|\big)~\leq~ \O(1)\cdot |y-x|\quad  \forall x,y\in\R^d.
\]
 In addition, there exists a constant  $r_0>0$ such that 
\bel{low-r}
r(x,\omega)~\geq~r_0\qquad\forall x\in \R^d, \omega\in U.
\eeq
\item [{\bf (H3).}] 
There exists a neighborhood $\mathcal{N}$ of $\partial\mathcal{S}$ such that 
\[
|g(y)-g(x)|~\leq~G\cdot |y-x|\quad\forall x,y\in \mathcal{N}, 
\]
with $0\leq G< \ds{r_0\over N}$.


\item [{\bf (H4).}] The target $\mathcal{S}$  has the $\rho_0$-inner ball property, i.e., for every $x\in\partial\mathcal{S}$, there exists $\xi_x\in N^{P}_{\overline{\mathcal{S}^c}}(x)$ with $|\xi_x|=1$ such that 
\[
\xi_x\cdot (y-x)~\leq~{1\over 2\rho_0}\cdot |y-x|^2\qquad\forall y\in \overline{\mathcal{S}^c}.
\] 

\end{itemize}
Under the assumption {\bf (H1)}, for any $x\in \R^d$ and $u\in\mathcal{U}_{ad}$, the system (\ref{Cont}) admits a unique Carath\'eodory  solution denoted by $y^{x,u}(\cdot)$. The exit time of $y^{x,u}(\cdot)$ is denoted by 
\[
\tau^{x,u}~\doteq~\min\big\{t\geq 0: y^{x,u}(t)\in \mathcal{S}\big\}
\]
with the convention that $\tau^{x,u}=+\infty$ if $y^{x,u}(t)\notin \mathcal{S}$ for all $t\geq 0$. The assumptions {\bf (H2)} and {\bf (H3)} imply that, given $x\in \R^d$ and $u\in\mathcal{U}_{ad}$, the function
\bel{si}\tau~\mapsto~\int_{0}^{\tau}r\left(y^{x,u}(s),u(s)\right)ds+g\left(y^{x,u}(\tau)\right)\eeq
is strictly increasing on a neighborhood of $\tau^{x,u}$. Without loss of generality, we shall suppose that  $g$ is Lipschitz on $\R^d$ with a Lipschitz constant $\ds G< {r_0\over N}$. The {\it controllable set} $\mathcal{R}$ is  the set of all points $x$ which can be steered to the target in a finite time i.e., 
\[
\mathcal{R}~=~\bigcup_{t>0}\mathcal{R}(t),\qquad\mathcal{R}(t)~=~\left\{x\in \R^d: T_{\mathcal{S}}(x)<t\right\},\qquad t>0,
\]
with $T_{\mathcal{S}}(\cdot)$ being the minimum time function
\[
T_{\mathcal{S}}(x)~\doteq~\inf_{u\in\mathcal{U}_{ad}}\tau^{x,u}\qquad\forall x\in\R^d.
\]
Notice that the value function $V$ defined in (\ref{V1}) is  the minimum time function $T_S$ in the case of $r\equiv 1, g\equiv 0$. From \cite[Section IV, Proposition 1.2]{BD}, if the system $(f,U)$ is small time controllable on $\mathcal{S}$ (briefly STC$\mathcal{S}$), i.e., 
\[
\mathcal{S}~\subseteq~\mathrm{int}(\mathcal{R}(t))\qquad\forall t>0,
\]
then $T_S$ is continuous on $\mathcal{R}$. In general, the assumption {\bf(H3)} can be regarded as a compatibility condition on the terminal cost $g$ which allows us to obtain the continuity of the value function $V$ together with the (STC$\mathcal{S}$) condition. From \cite[Section IV, Theorem~3.6 and Proposition~3.7]{BD}, one can easily derive the following lemma.
\begin{lemma}\label{Con-V} In addition to the standard assumptions {\bf (H1)-(H2)}, if both (STC$\mathcal{S}$) and {\bf (H3)} hold then $V$ is continuous in $\mathcal{R}$.
\end{lemma}

However,  if no restriction is made on the Lipschitz constant of $g$, the value function can be discontinuous along the boundary of the target $\mathcal{S}$. Indeed, the (STC$\mathcal{S}$) condition is not sufficient for the continuity of $V$  (see e.g. in \cite[Section IV, example 3.1 and example 3.2]{BD}).
Finally, the assumption {\bf (H0)} is related to the existence of optimal trajectories. In fact one can prove, by standard techniques, the following result.
\begin{theorem} Under hypotheses {\bf (H0)}-{\bf (H3)}, there exists a minimizer for the optimal control problem (\ref{Cont})-(\ref{V1}) for any choice of initial point $x\in\mathcal{R}$. Moreover, the uniform limit of optimal trajectories is an optimal trajectory; that is, if $y_n(\cdot)$
are trajectories converging uniformly to $y(\cdot)$ and every $y_n(\cdot)$ is optimal for the point
$x_n\doteq y_n(0)$, then $y(\cdot)$ is optimal for $x\doteq\ds\lim_{n\to\infty}x_n$.
\end{theorem}
The  Hamiltonian $H$ associated to (\ref{Cont})-(\ref{V1}) is defined by
\bel{H}
H(x,p)~\doteq~\max_{\omega\in U}\{ - p\cdot f(x,\omega)-r(x,\omega)\},\qquad (x,p)\in \R^d\times\R^d,
\eeq
while the horizontal Hamiltonian $H^0$ is defined by 
\bel{H0}
H^0(x,p)~\doteq~ \ds\max_{\omega\in U}\left\{-p\cdot f(x,\omega)\right\},\qquad (x,p)\in \R^d\times \R^d.
\eeq

We now prove a version of the Pontryagin maximum principle for  non-horizontal cases that generalizes the statement in \cite[Theorem 4.3]{PCC} for continuous value function and nonsmooth target. The horizontal cases will be studied in Section 4 by using the regularity results of Section 3.

\begin{lemma}\label{O}  In addition to {\bf (H1)}-{\bf (H3)}, assume that $g$ is locally semiconcave.  Let $x\in \mathcal{R}\backslash\mathcal{S}$ and let $u^*\in \mathcal{U}_{ad}$ be optimal for $x$. Set 
\[
y^*(\cdot)~\doteq~y^{x,u^*}(\cdot),\qquad \tau^*~\doteq~\tau^{x,u^*}\qquad\text{and}\qquad x^*~\doteq~y^*(\tau^*).
\]
Given $q^*\in D^+g(x^*)$, let $\xi^*$ be a proximal normal vector to $\overline{\mathcal{S}^c}$ at $x^*$ such that $H(x^*,p^*)=0$, with $p^*=q^*-\xi^*$. Let  $p:[0,\tau^*]\to\R^d$ be  the solution to the system
\bel{adjoint}
p'(t)~=~-D^{T}_xf(y^*(t),u^*(t)) \,p(t)-D_xr(y^*(t),u^*(t)),\qquad p(\tau^*)~=~p^*.
\eeq
Then for a.e. $t\in [0,\tau^*]$, it holds 
\bel{H-cond}
H(y^*(t),p(t))~=~- p(t)\cdot f(y^*(t),u^*(t))-r(y^*(t),u^*(t))~=~0.
\eeq
\end{lemma}
{\bf Proof.} {\bf 1.} We first show that the first inequality in (\ref{H-cond}) holds at every $t\in\, ]0,\tau^*[$ that is a Lebesgue point for the functions $f(y^*(\cdot),u^*(\cdot)), r(y^*(\cdot),u^*(\cdot))$. For a sake of simplicity, we shall assume that $t=0$, the computations in the general case being entirely analogous. Given $\omega\in U$, we show that 
\bel{Max1}
-p(0)\cdot f(x,\omega)-r(x,\omega)~\leq~-p(0)\cdot f(x,u^*(0))-r(x,u^*(0)).
\eeq Let $\omega^*\in U$ be such that 
\bel{H=00}
H(x^*,p^*)~=~-p^*\cdot f(x^*,\omega^*)-r(x^*,\omega^*)~=~0.
\eeq
Defining, for $\ve\in\,]0,\tau^*[$,
\[
u_{\ve}(s)~=~\left\{\begin{array}{llllll}\omega & s\in[0,\ve]\\u^*(s) & s\in\,]\ve,\tau^*[\\ \omega^* & s\in[\tau^*,+\infty[,\end{array}\right.
\]
we set $y_{\ve}(\cdot)\doteq y^{x,u_{\ve}}(\cdot)$ and $\tau_{\ve}\doteq \tau^{x,u_{\ve}}$. Since $0$ is a  Lebesgue point of $r(y^*(\cdot),u^*(\cdot))$ and $f(y^*(\cdot),u^*(\cdot))$,  one has
\bel{rLeb}
\int_{0}^{\ve}\big[r(y_{\ve}(s),\omega)-r(y^*(s),u^*(s))\big]ds~=~\ve\cdot \big[r(x,\omega)-r(x,u^*(0))\big]+o(\ve)
\eeq
and
\bel{}
y_{\ve}(s)~=~y^*(s)+\ve\cdot v(s)+o(\ve)\qquad\forall s\in [\ve,\tau^*],
\eeq
with $v(\cdot)$ being the solution to the linearized system
\bel{L-v}
v'(s)~=~D_xf(y^*(s),u^*(s))\,v(s),\qquad v(\ve)~=~f(x,\omega)-f(x,u^*(0)).
\eeq

Since $y^*(\cdot)$ is an optimal trajectory, it holds 
\bel{J-v-J}
\begin{split}
0&~\leq~\int_{0}^{\tau_{\ve}}r(y_{\ve}(s),u_{\ve}(s))ds-\int_{0}^{\tau^*}r(y^*(s),u^*(s))ds+g(y_{\ve}(\tau_{\ve}))-g(x^*).
\end{split}
\eeq
Two cases are considered: 
\medskip

{\bf Case 1.} Assume that $\tau_{\ve}\leq \tau^*$. By (\ref{rLeb}) and (\ref{J-v-J}), we first have 
\[
\begin{split}
r_0\cdot (\tau^*-&\tau_{\ve})~\leq~\int_{\tau_{\ve}}^{\tau^*}r(y^*(s),u^*(s))ds\\
&~\leq~\int_{0}^{\ve}\big[r(y_{\ve}(s),\omega)-r(y^*(s),u^*(s))\big]ds+\int^{\tau_\ve}_{\ve}\big[r(y_{\ve}(s),u^*(s))-r(y^*(s),u^*(s))\big]ds\\
&\qquad\qquad\qquad\qquad +\big|g(y_{\ve}(\tau_{\ve}))-g(y^*(\tau_{\ve}))\big|+ \big|g(y^*(\tau_{\ve}))-g(x^*)\big|\\
&~\leq~\O(1)\cdot\ve+NG (\tau^*-\tau_{\ve}),
\end{split}
\]
and this implies that  
\bel{tau-ve-1}
\tau^*-\tau_{\ve}~\leq~ {\O(1)\over r_0-GN}\cdot \ve.
\eeq
Consequently, it holds
\[
\begin{split}
0&~\leq~\ve\cdot \big[r(x,\omega)-r(x,u^*(0))\big]+\int^{\tau_\ve}_{\ve}\big[r(y_{\ve}(s),u^*(s))-r(y^*(s),u^*(s))\big]ds\\
&\qquad\qquad\qquad\qquad\qquad\qquad- \int_{\tau_{\ve}}^{\tau^*}r(y^*(s),u^*(s))ds+g(y_{\ve}(\tau_{\ve}))-g(x^*)+o(\ve)\\
&~=~\ve\cdot \big[r(x,\omega)-r(x,u^*(0))\big]+\int_{\ve}^{\tau_{\ve}}D_xr(y^*(s),u^*(s))\cdot \big(y_{\ve}(s)-y^*(s)\big)ds\\
&\qquad\qquad\qquad\qquad\qquad\qquad - \int_{\tau_{\ve}}^{\tau^*}r(x^*,u^*(s))ds+g(y_{\ve}(\tau_{\ve}))-g(x^*)+o(\ve).
\end{split}
\]
On the other hand, recalling that $H(x^*,p^*)=0$ and $\xi^*\in N^{P}_{\overline{\mathcal{S}^c}}(x^*)$, we have 
\[
\begin{split}
- \int_{\tau_{\ve}}^{\tau^*}&r(x^*,u^*(s))ds~\leq~p^*\cdot \int_{\tau_{\ve}}^{\tau^*}f(x^*,u^*(s))ds~\leq~p^*\cdot \big(x^*-y^*(\tau_{\ve})\big)+o(\ve)\\
&~\leq~p^*\cdot \big(y_{\ve}(\tau_{\ve}) -y^*(\tau_{\ve})\big)+p^*\cdot \big(x^*-y_{\ve}(\tau_{\ve})\big)+o(\ve)\\
&~\leq~p(\tau_{\ve})\cdot \big(y_{\ve}(\tau_{\ve}) -y^*(\tau_{\ve})\big)+(\xi^*-q^*)\cdot  \big(y_{\ve}(\tau_{\ve})-x^*\big)+o(\ve)\\
&~\leq~p(\ve)\cdot \big(y_{\ve}(\ve)-y^*(\ve)\big)+\int_{\ve}^{\tau_{\ve}}{d\over ds}\big[p(s)\cdot \big(y_{\ve}(s)-y^*(s)\big)\big]ds-q^*\cdot  \big(y_{\ve}(\tau_{\ve})-x^*\big)+o(\ve)\\
&~\leq~\ve\cdot p(0)\cdot v(\ve)+\int_{\ve}^{\tau_{\ve}}{d\over ds}\big[p(s)\cdot \big(y_{\ve}(s)-y^*(s)\big)\big]ds-q^*\cdot  \big(y_{\ve}(\tau_{\ve})-x^*\big)+o(\ve).
\end{split}
\]
Summing the above two estimates and using   $q^*\in D^+g(x^*)$, we obtain
\bel{ep}
\begin{split}
0&~\leq~\ve\cdot \big[r(x,\omega)+p(0)\cdot f(x,\omega)-r(x,u^*(0))-p(0)\cdot f(x,u^*(0))\big]\\
&~~~~ +\int_{\ve}^{\tau_{\ve}}\big[D_xr(y^*(s),u^*(s))+p'(s)\big]\cdot\big(y_{\ve}(s)-y^*(s)\big)ds\\
&~~~~ +\int_{\ve}^{\tau_{\ve}}p(s)\cdot \big[f(y_{\ve}(s),u^*(s))-f(y^*(s),u^*(s))\big]ds+o(\ve)\\
&~\leq~\ve\cdot \big[r(x,\omega)+p(0)\cdot f(x,\omega)-r(x,u^*(0))-p(0)\cdot f(x,u^*(0))\big]+o(\ve),
\end{split}
\eeq
and this yields (\ref{Max1}).
\medskip

{\bf Case 2.} Assume that $\tau_{\ve}> \tau^*$. In this case, if 
\[
\begin{split}
0&~\leq~\int_{0}^{\tau^*}r(y_{\ve}(s),u_{\ve}(s))ds-\int_{0}^{\tau^*}r(y^*(s),u^*(s))ds+g(y_{\ve}(\tau^*))-g(x^*),
\end{split}
\]
then one can follow the same argument in Case 1 to achieve (\ref{Max1}). Otherwise, there exists $\tau_1 \in\, ]\tau^*,\tau_{\ve}]$ such that 
\bel{tau-1}
0~=~\int_{0}^{\tau_1}r(y_{\ve}(s),u_{\ve}(s))ds-\int_{0}^{\tau^*}r(y^*(s),u^*(s))ds+g(y_{\ve}(\tau_1))-g(x^*).
\eeq
In this case, we estimate 
\[
\begin{split}
r_0\cdot (\tau_1-\tau^*)&~\leq~\int_{\tau^*}^{\tau_1}r(y_{\ve}(s),u_{\ve}(s))ds\\
&~=~\int_{0}^{\tau^*}\big[r(y^*(s),u^*(s))-r(y_{\ve}(s),u_{\ve}(s))\big]ds+g(y_{\ve}(\tau_1))-g(x^*)\\
&~\leq~\O(1)\cdot \ve+G\cdot |y_{\ve}(\tau_1)-x^*|~\leq~\O(1)\cdot\ve + GN\cdot (\tau_1-\tau^*)
\end{split}
\]
and this implies that 
\[
\tau_1-\tau^*~\leq~{\O(1)\over r_0-GN}\cdot \ve.
\]
Thus, from (\ref{tau-1}), we derive  
\[
\begin{split}
0&~=~\int_{0}^{\tau^*}\big[r(y_{\ve}(s),u_{\ve}(s))-r(y^*(s),u^*(s))\big]ds+\int_{\tau^*}^{\tau_{1}}r(y_{\ve}(s),u_{\ve}(s))ds+g(y_{\ve}(\tau_{1}))-g(x^*)\\
&~=~\int_{0}^{\ve}\big[r(y_{\ve}(s),\omega)-r(y^*(s),u^*(s))\big]ds+\int_{\ve}^{\tau^*}\big[r(y_{\ve}(s),u^*(s))-r(y^*(s),u^*(s))\big]ds\\
&\qquad\qquad\qquad\qquad\qquad +\int_{\tau^*}^{\tau_{1}}r(y_{\ve}(s),\omega^*)ds+q^*\cdot  \big(y_{\ve}(\tau_{1})-x^*\big)+o(\ve)\\
&~\leq~\ve\cdot \big[r(x,\omega)-r(x,u^*(0))\big]+
\int_{\ve}^{\tau^*}D_xr(y^*(s),u^*(s))\cdot\big(y_{\ve}(s)-y^*(s)\big)ds\\
&\qquad\qquad\qquad +(\tau_{1}-\tau^*)\cdot r(x^*,\omega^*)+(\tau_{1}-\tau^*)q^*\cdot f(x^*,\omega^*)+ q^*\cdot \big(y_{\ve}(\tau^*)-x^*\big)+o(\ve)
\end{split}
\]
%
In addition, since $\xi^*$ is a proximal normal vector to $\overline{\mathcal{S}^c}$ at $x^*$, one has
\[
\begin{split}
\xi^*\cdot \big(y_{\ve}(\tau^*)-x^*\big)&~=~\xi^*\cdot \big(y_{\ve}(\tau^*)-y_{\ve}(\tau_{1})\big)+\xi^*\cdot\big(y_{\ve}(\tau_{1})-x^*\big)+o(\ve)\\
&~\leq~-(\tau_{1}-\tau^*)\xi^*\cdot f(x^*,u^*)+\O(1)\cdot |y_{\ve}(\tau_{1})-x^*|^2+o(\ve)\\
&~=~-(\tau_{1}-\tau^*)\xi^*\cdot f(x^*,u^*)+o(\ve),
\end{split}
\]
and (\ref{H=00}) yields 
\[
\begin{split}
0&~\leq~\ve\cdot \big[r(x,\omega)-r(x,u^*(0))\big]+\int_{\ve}^{\tau^*}D_xr(y^*(s),u^*(s))\cdot  \big(y_{\ve}(s)-y^*(s)\big)ds\\
&~~~~\qquad\qquad \qquad\qquad\qquad+p(\tau^*) \cdot  \big(y_{\ve}(\tau^*)-y^*(\tau^*)\big)+o(\ve).
\end{split}
\]
Thus, one can follow the same estimates in (\ref{ep}) to derive (\ref{Max1}).
\v

{\bf 2.} Finally, we observe that for every $s,t\in  [0, \tau^*]$, with $t$ satisfying the first inequality in (\ref{H-cond}), it holds 
\[
\begin{split}
H(y^*(t),p(t))&- H(y^*(s),p(s))~\leq~\big(p(s)-p(t)\big) \cdot f(y^*(t),u^*(t))\\
&~~~~~+p(s)\cdot \big[f(y^*(s),u^*(t))-f(y^*(t),u^*(t))\big]+r(y^*(s),u^*(t))-r(y^*(t),u^*(t))\\
&~=~\big(p(s)-p(t)\big) \cdot f(y^*(t),u^*(t))-p'(t)\cdot \big(y^*(s)-y^*(t)\big)+o(|t-s|).
\end{split}
\]
 Moreover, if  $p(\cdot),y^*(\cdot)$ and $ H(y^*(\cdot),p(\cdot))$ are differentiable at $t$ then $\ds{d\over dt}H(y^*(t),p(t))=0$. Indeed, in this case, one has
\[
\begin{split}
{d\over dt}H(y^*(t),p(t))&~=~\lim_{s\to t+} {H(y^*(s),p(s))-H(y^*(t),p(t))\over s-t}\\
&~\geq~\lim_{s\to t+}\Big[{p(t)-p(s)\over s-t} \cdot f(y^*(t),u^*(t))+p'(t)\cdot {y^*(s)-y^*(t)\over s-t}+o(1)\Big]\\
&~=~0
\end{split}
\]
and 
\[
\begin{split}
{d\over dt}H(y^*(t),p(t))&~=~\lim_{s\to t-} {H(y^*(s),p(s))-H(y^*(t),p(t))\over s-t}\\
&~\leq~\lim_{s\to t-}\Big[{p(t)-p(s)\over s-t} \cdot f(y^*(t),u^*(t))+p'(t)\cdot {y^*(s)-y^*(t)\over s-t}+o(1)\Big]\\
&~=~0.
\end{split}
\]
Finally, since the map $t\mapsto H(y^*(t),p(t))$ is absolutely continuous and $H(y^*(\tau^*),p(\tau^*))=0$,  one obtains $H(y^*(t),p(t))=0$ for all $t\in [0,\tau^*]$.
\endproof
\v

To complete this subsection, let us define the reachable boundary of target $\mathcal{S}$ as
\bel{S*}
\begin{split}
\partial\mathcal{S}^{*}~\doteq~\{x^*\in\partial\mathcal{S}&: \exists (x^*_n)_{n\geq 1} \text{ in }\partial\mathcal{S}, \,\exists\delta>0~\text{such that}~\lim_{n\to\infty}x^*_n=x^*~\text{and}\\
& x^*_n~\text{is~reached~by~a~trajectory~of~}(\ref{Cont})~\text{in~an~amount~of~time}~\geq \delta \}.
\end{split}
\eeq
For every $x^*\in\partial\mathcal{S}^*$,  we shall denote
\bel{Nor}
\begin{split}
N_0(x^*)&~\doteq~\left\{0\neq \xi\in N^{P}_{\overline{\mathcal{S}^c}}(x^*): \max_{\omega\in U}~\xi\cdot f(x^*,\omega)=0\right\},\\
N_{1}(x^*)&~\doteq~\left\{\xi\in N^{P}_{\overline{\mathcal{S}^c}}(x^*): \max_{\omega\in U}~\xi\cdot f(x^*,\omega)>0\right\}.
\end{split}
\eeq
Notice that $N_0(x^*)$ and $ N_{1}(x^*)$ can be both nonempty. The Petrov condition (\ref{Petrov}) is satisfied at $x^*$ if and only if  the set $ N_{0}(x^*)$ is empty. Moreover, the following property holds:
\begin{lemma}\label{NN} For every $x^*\in\partial\mathcal{S}^*$, it holds 
\bel{NN1}
N_0(x^*)\cup  N_{1}(x^*)~=~ N^{P}_{\overline{\mathcal{S}^c}}(x^*)\smallsetminus\{0\}.
\eeq
\end{lemma}
{\bf Proof.} Fix  $x^*\in \partial\mathcal{S}^*$. By the definition of $\partial\mathcal{S}^*$, there exist a constant $\delta>0$ and a sequence of trajectories $y^{x_n,u_n}(\cdot)$ of (\ref{Cont}) with $x_n^*\doteq y^{x_n,u_n}(\tau^{x_n,u_n})$ such that 
\[
\lim_{n\to\infty} x^*_n~=~x^*,\qquad  \delta~\leq~ \tau^{x_n,u_n}<+\infty\qquad\forall n\geq 1.
\]
For every  $ 0\neq \xi\in N^{P}_{\overline{\mathcal{S}^c}}(x^*)$, it holds
\begin{multline*}
\xi\cdot \big(y^{x_n,u_n}~(\tau^{x_n,u_n}-\ve)-x^*\big)~\leq~\O(1)\cdot \left|y^{x_n,u_n}(\tau^{x_n,u_n}-\ve)-x^*\right|^2\\
~\leq~\O(1)\cdot \left(\big|y^{x_n,u_n}(\tau^{x_n,u_n}-\ve)-x_n^*\big|^2+\big| x_n^*-x^*\big|^2\right)~\leq~\O(1)\cdot \left(\ve^2+\big| x_n^*-x^*\big|^2\right)
\end{multline*}
for all $n\geq 1$, $\ve\in \,]0,\delta[$. On the other hand,
\[
\begin{split}
\xi\cdot \left(y^{x_n,u_n}(\tau^{x_n,u_n}-\ve)-x^*\right)&~\geq~-\int_{\tau^{x_n,u_n}-\ve}^{\tau^{x_n,u_n}}\xi\cdot f(x^*,u_n(s))ds-\O(1)\cdot (|x^*_n-x^*|+\ve^2)\\
&~\geq~-\ve\cdot \max_{\omega\in U} \,\xi\cdot f(x^*,\omega)-\O(1)\cdot (|x^*_n-x^*|+\ve^2).
\end{split}
\]
We then derive 
\[
\max_{\omega\in U}~\xi\cdot f(x^*,\omega)~\geq~-\O(1)\cdot  \left({|x^*_n-x^*|\over \ve}+\ve\right).
\]
Thus, choosing $\ve= |x^*_n-x^*|^{1/2}$ and taking $n\to +\infty$, we obtain 
$$\max_{\omega\in U}~\xi\cdot f(x^*,\omega)~\geq~0,$$
and this yields (\ref{NN1}).
\endproof
\medskip

Finally, let us introduce the sets of {\it transported proximal normals}  $N^{\dagger}_1(x)$ and  {\it transported horizontal proximal normals}  $N^{\dagger}_0(x)$ at any point  $x\in  \mathcal{S}^c$ according with  $N_0(x^*)$ and $N_1(x^*)$. These sets play an important role in the representation formula of proximal  supergradients and proximal horizontal supergradients of the value function $V$.
\begin{definition}\label{TP}
For every $x\in \mathcal{S}^c$, calling
\bel{Ux}
\mathcal{U}_x~\doteq~\left\{u^*\in\mathcal{U}_{ad}:u^*~\text{is an optimal control which steers $x$ to $\mathcal{S}$}\right\},
\eeq
we shall  denote by 
\begin{itemize}
\item  [$\diamond$] $N^{\dagger}_0(x)$ the set of  {\it transported horizontal proximal normals}, i.e., the set of $p(0)\in \R^d$ such that $p(\cdot)$ is the unique solution of the ODE
\bel{adjoint-22}
p'(t)~=~-D^{T}_xf(y^{x,u^*}(t),u^*(t))\cdot p(t),\qquad p(\tau^{x,u^*})~=~-\xi^*,
\eeq
for some  $u^*\in\mathcal{U}_x$ and $\xi^*\in N_0\big(y^{x,u^*}(\tau^{x,u^*})\big)$;

\item  [$\diamond$]  $N^{\dagger}_1(x)$ the set of {\it transported  proximal normals}, i.e., the set of $p(0)\in \R^d$ such that $p(\cdot)$ is the unique solution of the ODE
\bel{adjoint-12}
p'(t)~=~-D^{T}_xf(y^{x,u^*}(t),u^*)\cdot p(t)-D_xr(y^{x,u^*}(t),u^*(t)),\qquad p(\tau^{x,u^*})~=~p^*,
\eeq 
for some  $u^*\in\mathcal{U}_x$, $\xi^*\in N_1(y^{x,u^*}(\tau^{x,u^*}))$ and  $p^*=q^*-\xi^*$ such that 
$$
H(y^{x,u^*}(\tau^{x,u^*}), p^*)~=~0,\qquad q^*\in D^+g(y^{x,u^*}(\tau^{x,u^*})).
$$
\end{itemize}

\end{definition}
\section{Regularity theory for the value function}
\setcounter{equation}{0}
In this section, we study the regularity of the value function $V$ under a weak controllability condition such that the value function $V$ is just continuous. From the viewpoint of Proposition \ref{hypo-upp}, we first show that  the hypograph of $V$ satisfies an exterior sphere condition.

\begin{theorem}\label{Main1} Assume that  {\bf (H0)}-{\bf (H4)} hold and $g$ is locally semiconcave. If the value function $V$ is continuous then $\mathrm{hypo}(V)$ satisfies an exterior sphere condition and 
\[
N_1^{\dagger}(x)\bigcup N_0^{\dagger}(x)~\neq~\varnothing\qquad\forall x\in \mathcal{R}\backslash\mathcal{S}.
\]
Moreover, the proximal (horizontal) superdifferential of  $V$ at $x$ satisfies
\bel{sup-V}
\partial^{P}V(x)~\supseteq~ \partial^{\infty}V(x)+\mathrm{co}\big[N^{\dagger}_1(x)\big],\qquad \partial^{\infty}V(x)~\supseteq~\mathrm{co}\big[ N^{\dagger}_0(x)\cup \{0\}\big].
\eeq
\end{theorem}
As a consequence, $V$ satisfies all the properties listed in Proposition \ref{hypo-upp}. Moreover, from \cite{ANV,N}, it holds:
\begin{corollary}\label{semic}
Under the same assumption in Theorem \ref{Main1}, if   the Petrov condition (\ref{Petrov}) holds then $V$ is locally semiconcave.
\end{corollary}

In addition, from Proposition \ref{EQ}, one obtains the following result which extends  Theorem 3.3 in \cite{GN} to the more general cases of  optimal exit time problems. 

\begin{corollary} Under the same assumption in Theorem \ref{Main1}, if $N^P_{\mathrm{hypo}(V)}(x,V(x))$ is pointed, i.e., 
\bel{pointed}
N^P_{\mathrm{hypo}(V)}(x,V(x))\cap \left(-N^P_{\mathrm{hypo}(V)}(x,V(x))\right)~=~\{0\}.
\eeq
  for every $x\in \mathcal{R}\backslash\mathcal{S}$ then $\mathrm{hypo}(V)$ is a $\vp$-convex set. 
\end{corollary}
An example  showing that the above Corollary  is no longer valid if the pointedness assumption is dropped in the case of the minimum time function for nonlinear control systems is provided in \cite[Example 2]{GN}. In this case, the example also points out that the first inclusion in (\ref{sup-V}) can be strict in general. 
\medskip

Next we shall establish a representation formula of  supergradients of $V$ and proximal horizontal supergradients of $V$ at $x$ by using the sets of transported normals  under the pointedness assumption on  $N^P_{\mathrm{hypo}(V)}(x,V(x))$ under different conditions.
%
%
%
Our second main results  are stated as follows:
\begin{theorem}\label{Re-pre} Under the same hypotheses of Theorem \ref{Main1}, let $x\in\mathrm{int}\mathcal{R}\backslash\mathcal{S}$.
\begin{itemize}
\item [(i).] Assume that $H$ is $C^{1}$ on  $\R^d\times \big(\R^d\backslash\{0\}\big)$ then $\partial^{*}V(x)\subseteq N_1^{\dagger}(x)$. If in addition   $N^P_{\mathrm{hypo}(V)}(x,V(x))$ is pointed, it holds
 \bel{sup-V-e}
\partial^{P}V(x)~=~ \partial^{\infty}V(x)+\mathrm{co}\big[N^{\dagger}_1(x)\big].
\eeq

\item [(ii).] Assume that $H^0$ is $C^{1}$ on  $\R^d\times \big(\R^d\backslash\{0\}\big)$ then $\partial^{*,\infty}V(x)\subseteq N_0^{\dagger}(x)$.  
If in addition $N^P_{\mathrm{hypo}(V)}(x,V(x))$ is pointed, it holds
 \bel{sup-V-ii}
\partial^{\infty}V(x)~=~\mathrm{co}\big[ N^{\dagger}_0(x)\cup \{0\}\big].
\eeq
\end{itemize}
\end{theorem}

\begin{corollary} Under the same setting in Theorem \ref{Main1}, if both $H$ and $H^0$ are $C^{1}$ on  $\R^d\times \big(\R^d\backslash\{0\}\big)$ and $x\in\mathrm{int}\mathcal{R}\backslash\mathcal{S}$ is such that $N^P_{\mathrm{hypo}(V)}(x,V(x))$ is pointed, then 
\[
\partial^{P}V(x)~=~\mathrm{co}\big[N^{\dagger}_1(x)\big]+\mathrm{co}\big[ N^{\dagger}_0(x)\cup \{0\}\big],\qquad \partial^{\infty}V(x)~=~\mathrm{co}\big[ N^{\dagger}_0(x)\cup \{0\}\big].
\]
\end{corollary}

As a consequence of Corollary \ref{semic} and Theorem \ref{Re-pre},  since proximal normals to the hypograph of a semiconcave function is pointed at every point, one obtains the following result  for $V$ under the Petrov condition. 

\begin{corollary}
Under the same assumptions in Theorem \ref{Main1}, if  $H$ is $C^{1}$ on  $\R^d\times \big(\R^d\backslash\{0\}\big)$ and the Petrov condition (\ref{Petrov}) holds then  
\[
D^+V(x)~=~\partial V^P(x)~=~\mathrm{co}\big[N^{\dagger}_1(x)\big].
\]
for all $x\in  \mathrm{int}\mathcal{R}\backslash\mathcal{S}$.
\end{corollary}
\begin{remark} In general we cannot expect the maps $H(x,p)$ and $H^0(x,p)$ to be differentiable when $p=0$. On the other hand, their differentiability on $\R^d\times \big(\R^d\backslash\{0\}\big)$ can be deduced by \cite[Theorem 3.4.4]{PC} in some special cases. In particular, if  for all $(x,p)\in \R^d\times \big(\R^d\backslash\{0\}\big)$, there exists a unique $u(x,p)\in U$ such that 
\[
u(x,p)~\in~\mathrm{argmax}_{\omega\in  U} [-f(x,\omega)\cdot p-r(x,\omega)],
\]
and the map $(x,p)\mapsto u(x,p)$ is continuous in $\R^d\times \big(\R^d\backslash\{0\}\big)$, then $H$ is $C^{1}$ on  $\R^d\times \big(\R^d\backslash\{0\}\big)$.
\\
Similarly, if  for all $(x,p)\in \R^d\times \big(\R^d\backslash\{0\}\big)$, there exists a unique $u_0(x,p)\in U$ such that 
\[
u_0(x,p)~\in~\mathrm{argmax}_{\omega\in  U} [-f(x,\omega)\cdot p],
\]
and the map $(x,p)\mapsto u_0(x,p)$ is continuous in $\R^d\times \big(\R^d\backslash\{0\}\big)$, then $H^0$ is $C^{1}$ on  $\R^d\times \big(\R^d\backslash\{0\}\big)$.
\end{remark}
%
\subsection{Transported proximal normal vectors }
Toward a proof of Theorem \ref{Main1} and Theorem \ref{Re-pre}, we shall establish several results on transported proximal normal vectors defined in (\ref{adjoint-22})-(\ref{adjoint-12}). In particular, Proposition \ref{Pro-0} and Proposition \ref{Pro-1} yield (\ref{sup-V}). The first result is concerned with optimal trajectories which reach the target $\mathcal{S}$ tangentially.
\begin{proposition}\label{Pro-0}  In addition to {\bf (H1)}-{\bf (H3)},  suppose that $V$ is upper  semicontinuous. Given $x\in\mathcal{S}^c$, let $u^*(\cdot)$ be an optimal control which steers $x$ to $x^*=y^{x,u^*}(\tau^*)\in\partial\mathcal{S}$ with $\tau^*\doteq\tau^{x,u^*}$ and let $\xi^*\in N^{P}_{\overline{\mathcal{S}^c}}(x^*)$ be realized by a ball of radius $\bar{\rho}>0$. If $\xi^*\in N_0(x^*)$ then, calling $p(\cdot)$  the solution of the ODE
\bel{adjoint-23}
p'(t)~=~-D^{T}_xf(y^{x,u^*}(t),u^*(t))\cdot p(t),~\qquad p(\tau^*)~=~-\xi^*,
\eeq
one has  
\bel{p-0-V}
p(t)~\in~\partial^{\infty} V(y^{x,u^*}(t)),
\eeq
with 
\bel{hyp-f-0}
{-p(t)\over |p(t)|}\cdot \big(z-y^{x,u^*}(t)\big)~\leq~\O(1)\cdot \left(1+{1\over 2\bar{\rho}}\right)\cdot\left(|z-y^{x,u^*}(t)|^2+|\beta-V(y^{x,u^*}(t))|^2\right)
\eeq
for all $t\in [0,\tau^*]$ and $(z,\beta)\in \mathrm{hypo}(V)$.
\end{proposition}
{\bf Proof.} For simplicity, we prove the statement for $t=0$. We have 
\bel{V-x1}
\alpha~\doteq~V(x)~=~\int_{0}^{\tau^*}r\left(y^{x,u^*}(s), u^*(s)\right)ds+g(x^*).
\eeq
Given $ z\in\mathcal{S}^c$, set $\tau_z\doteq\tau^{z,u^*}$. We consider two cases:
\medskip

{\bf Case 1.} If $\tau_z> \tau^*$ then from (\ref{adjoint-23}) and {\bf (H1)}, it holds
\[
\begin{split}
\int^{\tau^*}_{0}{d\over ds}\left[ p(s)\cdot \big( y^{z,u^*}(s)-y^{x,u^*}(s)\big)\right]ds&~\leq~\O(1)\cdot\int_{0}^{\tau^*}|p(s)|\cdot \big|y^{z,u^*}(s)-y^{x,u^*}(s)\big|^2ds\\
&~\leq~\O(1)\cdot |p(0)|\cdot |x-z|^2.
\end{split}
\]
Setting $z_1\doteq y^{z,u^*}(\tau^*)\in\mathcal{S}^c$, we estimate 
\[
\begin{split}
-p(0)\cdot (&z-x)~=~-p(\tau^*)\cdot (z_1-x^*)+\int^{\tau^*}_{0}{d\over ds}\left[ p(s)\cdot \big( y^{z,u^*}(s)-y^{x,u^*}(s)\big)\right]ds\\
&\leq~{ |\xi^*|\over 2\bar{\rho}}\cdot |z_1-x^*|^2+\O(1)\cdot |p(0)|\cdot |x-z|^2~\leq~\O(1)\cdot |p(0)|\cdot \left(1+{1\over 2\bar{\rho}}\right)\cdot |x-z|^2.
\end{split}
\]
This yields (\ref{hyp-f-0}).

{\bf Case 2.}  If $\tau_z\leq \tau^*$ then $z^*\doteq y^{z,u^*}(\tau_z)~\in~\partial\mathcal{S}$ and, taking $\beta\leq V(z)$, we have
\[
 \beta~\leq~V(z)~\leq~\int_{0}^{\tau_z}r\left(y^{z,u^*}(s),u^*(s)\right)ds+g(z^*).
\]
From  {\bf (H2)} and (\ref{V-x1}), one has  
\[
\begin{split}
\int_{\tau_z}^{\tau^*}r(y^{x,u^*}(s),u^*(s))ds&~\leq~\int_{0}^{\tau_z}\big[r(y^{z,u^*}(s),u^*(s))-r(y^{x,u^*}(s),u^*(s))\big]ds+ g(z^*)-g(x^*)+\alpha-\beta\\
&~\leq~\O(1)\cdot |z-x|+GN\cdot(\tau^*-\tau_z)+\alpha-\beta,
\end{split}
\]
and {\bf (H3)} implies that 
\[
\tau^*-\tau_z~\leq~{1\over r_0-GN}\cdot \big[\O(1)\cdot |z-x|+\alpha-\beta\big].
\]
In particular, setting $x_1\doteq y^{x,u^*}(\tau_z)$, we  have 
\bel{T-est0}
|x^*-x_1|~\leq~ N\cdot (\tau^*-\tau_z)~\leq~{N\over r_0-GN}\cdot  \big[\O(1)\cdot |z-x|+\alpha-\beta\big]\,.
\eeq
Again from  (\ref{adjoint-23})  and {\bf (H1)} we have
\[
\begin{split}
\int^{\tau_z}_{0}{d\over ds}\left[ p(s)\cdot \big( y^{z,u^*}(s)-y^{x,u^*}(s)\big)\right]ds&~\leq~\O(1)\cdot\int_{0}^{\tau_z}|p(s)|\cdot \big|y^{z,u^*}(s)-y^{x,u^*}(s)\big|^2ds\\
&~\leq~\O(1)\cdot |p(0)|\cdot |x-z|^2.
\end{split}
\]
Then we compute
\begin{multline}\label{et0}
-p(0)\cdot (z-x)~=~-p(\tau_z)\cdot (z^*-x_1)+\int^{\tau_z}_{0}{d\over ds}\left[ p(s)\cdot \big( y^{z,u^*}(s)-y^{x,u^*}(s)\big)\right]ds\\
~\leq~ -p(\tau^*)\cdot (z^*-x_1)+|p(\tau^*)-p(\tau_z)|\cdot |z^*- x_1|+\O(1)\cdot |p(0)|\cdot |x-z|^2\\
~\leq~\xi^*\cdot (z^*-x^*)+\xi^*\cdot (x^*-x_1)+\O(1)\cdot |p(0)|\cdot \big( |z-x|^2+(\alpha-\beta)^2\big)\\
~\leq~\O(1)\cdot \left(1+{1\over 2\bar{\rho}}\right)\cdot |p(0)|\cdot \big( |z-x|^2+(\alpha-\beta)^2\big)+\int_{\tau_z}^{\tau^*}\xi^*\cdot f(y^{x,u^*}(s),u^*(s))ds.
\end{multline}
Since $\xi^*\in N_0(x^*)$,  for every $s\in [\tau_z,\tau^*]$ it holds
\[
\begin{split}
\xi^*\cdot f(y^{x,u^*}(s),u^*(s))&~\leq~\O(1)\cdot |\xi^*|\cdot |y^{x,u^*}(s)-x^*|+\max_{u\in U}\big\{\xi^*\cdot f(x^*,u)\big\}\\
&~\leq~\O(1)\cdot |\xi^*|\cdot (\tau^*-\tau_z)~\leq~\O(1)\cdot |p(0)|\cdot \big(\O(1)\cdot |z-x|+\alpha-\beta\big).
\end{split}
\]
Combining (\ref{T-est0}) and (\ref{et0}), we obtain
\[
-p(0)\cdot (z-x)~\leq~\O(1)\cdot \left(1+{1\over 2\bar{\rho}}\right)\cdot |p(0)|\cdot \big( |z-x|^2+(\alpha-\beta)^2\big),
\]
and this yields (\ref{hyp-f-0}).
\endproof

The rest of the subsection concerns optimal trajectories that reach the target $\mathcal{S}$ transversally.
 For every $\alpha>0$,  we show that the $\alpha$-sublevel of $V$,  denoted by 
\bel{Sr}
\mathcal{V}(\alpha)~=~\left\{x\in\R^d: V(x)\leq \alpha\right\}\,
\eeq
 satisfies an interior sphere condition under the same assumptions in Lemma \ref{O}. Such regularity property is propagated along the dual arc from the semiconcavity of $g$ and the internal sphere condition of $\mathcal{S}$. This will be achieved by the following two lemmas. 

\begin{lemma}\label{p1}   In addition to {\bf (H1)}-{\bf (H3)},  assume that $g$ is locally semiconcave and $V$ is upper semicontinuous. Given $x\in\mathcal{S}^c$, let $u^*(\cdot)$ be an optimal control which steers $x$ to $x^*=y^{x,u^*}(\tau^*)\in\partial\mathcal{S}$ with $\tau^*\doteq\tau^{x,u^*}$ and let $\xi^*\in N^{P}_{\overline{\mathcal{S}^c}}(x^*)$ be realized by a ball of radius $\bar{\rho}>0$. If $\xi^*\in N_1(x^*)$ and $q^*\in  D^+g(x^*) $ is such that  $H(x^*,p^*)=0$ with $p^*=q^*-\xi^*$ then, calling $p(\cdot):[0,\tau^*]\to\R^d$   the solution to the system
\bel{adjoint-1}
p'(t)~=~-D^{T}_xf(y^{x,u^*}(t),u^*(t))\cdot p(t)-D_xr(y^{x,u^*}(t),u^*(t)),\qquad p(\tau^*)~=~p^*,
\eeq 
one has $-p(0)\in N^{P}_{\overline{\mathcal{V}(\alpha)^c}}(x)$ with $\alpha\doteq V(x)$ and 
\bel{P-0}
{-{p(0)\over |p(0)|}}\cdot (z-x)~\leq~\O(1)\cdot \left(1+{1\over 2\bar{\rho}}\right)\cdot |z-x|^2\qquad\forall z\in {\overline{\mathcal{V}(\alpha)^c}}.
\eeq
\end{lemma}
{\bf Proof.} 
Since $u^*(\cdot)$ is an optimal control we have
\bel{V-x-1}
\alpha~=~V(x)~=~\int_{0}^{\tau^*}r\left(y^{x,u^*}(s), u^*(s)\right)ds+g(x^*).
\eeq
and by Lemma \ref{O}, for a.e. $t\in [0,\tau^*]$ it holds
\bel{H=0-1}
H(y^{x,u^*}(t),p(t))~=~- p(t)\cdot f(y^{x,u^*}(t),u^*(t))-r(y^{x,u^*}(t),u^*(t))~=~0.
\eeq
Since $\xi^*\in N^{P}_{\overline{\mathcal{S}^c}}(x^*)$ is realized by a ball of radius $\bar{\rho}>0$ and $q^*\in  D^+g(x^*) $, we  have 
\bel{inner-ball1}
\begin{cases}
\xi^*\cdot (y- x^*)~\leq~\ds{|\xi^*|\over 2\bar{\rho}}\cdot |y-x^*|^2\\[4mm]
g(y)-g(x^*)-q^*\cdot (y- x^*)~\leq~\O(1)\cdot |y-x^*|^2
\end{cases}
\qquad\forall y\in\overline{\mathcal{S}^c}.
\eeq
For any given $z\in  {\overline{\mathcal{V}(\alpha)^c}}$, set $ \tau_z\doteq\tau^{z,u^*}$. Two cases are considered:
\medskip

{\bf Case 1.} If $\tau_z\leq \tau^*$, set $z^*\doteq y^{z,u^*}(\tau_z)\in\partial\mathcal{S}$ and $x_1\doteq y^{x,u^*}(\tau_z)$.  Since
\[
 \alpha~\leq~V(z)~\leq~\int_{0}^{\tau_z}r\left(y^{z,u^*}(s),u^*(s)\right)ds+g(z^*),
\]
we obtain
\[
\begin{split}
r_0\cdot (\tau^*-&\tau_z)~\leq~\int_{\tau_z}^{\tau^*}r(y^{x,u^*}(s),u^*(s))ds\\
&~\leq~\int_{0}^{\tau_z}\big[r(y^{z,u^*}(s),u^*(s))-r(y^{x,u^*}(s),u^*(s))\big]ds+g(z^*)-g(x^*)\\
&~\leq~\O(1)\cdot|z-x|+NG \cdot(\tau^*-\tau_z),
\end{split}
\]
This implies 
\bel{tau-ve-12}
\tau^*-\tau_z~\leq~ {\O(1)\over r_0-GN}\cdot |z-x|
\eeq
and
\bel{T-est1}
|x^*-x_1|~\leq~ N\cdot (\tau^*-\tau_z)~\leq~{\O(1)\over r_0-GN}\cdot |z-x|\,.
\eeq
Hence,  we estimate 
\[
\begin{split}
-p(0)\cdot (z-x)&~=~-p(\tau_z)\cdot (z^*-x_1)+\int^{\tau_z}_{0}{d\over ds}\left[ p(s)\cdot \big( y^{z,u^*}(s)-y^{x,u^*}(s)\big)\right]ds\\
&~=~\big(p(\tau^*)-p(\tau_{z})\big)\cdot (z^*-x_1) -(q^*-\xi^*)\cdot (z^*-x^*)\\
&\qquad\qquad~~-p(\tau^*)\cdot (x^*-x_1) +\int^{\tau_z}_{0}{d\over ds}\left[ p(s)\cdot \big( y^{z,u^*}(s)-y^{x,u^*}(s)\big)\right]ds\\
&~\leq~\O(1)\cdot \left(1+{|\xi^*|\over 2\bar{\rho}}\right)\cdot |z^*-x^*|^2-q^*\cdot (z^*-x^*)\\
&\qquad\qquad~~-p(\tau^*)\cdot (x^*-x_1) +\int^{\tau_z}_{0}{d\over ds}\left[ p(s)\cdot \big( y^{z,u^*}(s)-y^{x,u^*}(s)\big)\right]ds\,.
\end{split}
\]
Recalling (\ref{V-x-1}), we compute 
\[
\begin{split}
\int^{\tau_z}_{0}{d\over ds}&\left[ p(s)\cdot \big( y^{z,u^*}(s)-y^{x,u^*}(s)\big)\right]ds\\
&~\leq~\O(1)\cdot |p^*|\cdot |z-x|^2-\int_{0}^{\tau_z}\ D_xr(y^{x,u^*}(s),u^*(s))\cdot \big(y^{z,u^*}(s)-y^{x,u^*}(s)\big) ds\\
&~\leq~\O(1)\cdot (1+|p^*|)\cdot |z-x|^2+\int_{0}^{\tau_z} r(y^{x,u^*}(s),u^*(s))ds-\int_{0}^{\tau_z} r(y^{z,u^*}(s),u^*(s))ds\\
&~\leq~\O(1)\cdot (1+|p^*|)\cdot |z-x|^2-\int_{\tau_z}^{\tau^*}r(y^{x,u^*}(s),u^*(s))ds+ g(z^*)-g(x^*),
\end{split}
\]
and  (\ref{H=0-1}) implies   
\[
\begin{split}
\int^{\tau_z}_{0}{d\over ds}&\left[ p(s)\cdot \big( y^{z,u^*}(s)-y^{x,u^*}(s)\big)\right]ds\\
&~\leq~\O(1)\cdot (1+|p^*|)\cdot |z-x|^2+\int_{\tau_z}^{\tau^*} p(s)\cdot f(y^{x,u^*}(s),u^*(s)) ds+g(z^*)-g(x^*).
\end{split}
\]
Noting that $H(x^*,p^*)=0$ implies $|p^*|\geq\ds{r_0\over N}$, we obtain 
\[
-p(0)\cdot (z-x)~\leq~\O(1)\cdot \left(1+{1\over 2\bar{\rho}}\right)\cdot |p^*| \cdot |z-x|^2+g(z^*)-g(x^*)-q^*\cdot (z^*-x^*),
\]
 and the seminconcavity property of  $g$ yields
\bel{in-est1}
-p(0)\cdot (z-x)~\leq~\O(1)\cdot \left(1+{1\over 2\bar{\rho}}\right)\cdot |p(0)|\cdot |z-x|^2.
\eeq

{\bf Case 2.} Assume that  $\tau_z> \tau^*$.  Let $\omega^*\in U$ be  such that 
\bel{chos1}
H(x^*,p^*)~=~-p^*\cdot f(x^*,\omega^*)-r(x^*,\omega^*)~=~0\,.
\eeq
We extend the optimal control $u^*$ to the interval $ ]\tau^*,+\infty[$ by setting 
\[
u^*(s)~=~\omega^*\qquad\forall s\in \,]\tau^*,+\infty[\,.
\]
Set $\tau^*_z\doteq\tau^{z,u^*}$ and consider the map $h_z:[0,{\tau^*_z}]\to\R$ defined by
\[
h_z(t)~=~\int_{0}^{t}r(y^{z,u^*}(s),u^*(s))ds+g(y^{z,u^*}(t))\qquad\forall t\in [0,{\tau^*_z}].
\]
By the condition {\bf (H3)},  for any $0<t_1<t_2<\tau^*_z$ we compute 
\[
h_z(t_2)-h_z(t_1)~=~\int_{t_1}^{t_2}r(y^{z,u^*}(s),u^*(s))ds+g(y^{z,u^*}(t_2))-g(y^{z,u^*}(t_1))~\geq~(r_0-GN)(t_2-t_1)~>~0,
\]
Therefore,  the map $t\mapsto h_z(t)$ is strictly increasing and 
\[
\alpha~=~V(x)~\leq~V(z)~\leq~h_z(\tau^*_z)\,.
\] 
Set $z_1\doteq y^{z,u^*}(\tau^*)$. By the same computations as in Case 1, we have 
\bel{et-11}
\begin{split}
- p(&0)\cdot (z-x)~=~-p^*\cdot (z_1-x^*)+\int^{\tau^*}_{0}{d\over ds}\Big[ p(s)\cdot\left(y^{z,u^*}(s)-y^{x,u^*}(s)\right)\Big] ds\\
&~\leq~ \O(1)\cdot |p(0)|\cdot |z-x|^2+ (\xi^*-q^*)\cdot (z_1-x^*)+g(z_1)-g(x^*)+\alpha-h_z(\tau^*).
\end{split}
\eeq
Two sub-cases are considered: 
\begin{itemize}
\item If $ h_z(\tau^*)\geq \alpha$ then 
\[
\begin{split}
- p(0)\cdot (z-x)&~\leq~\O(1)\cdot |p(0)|\cdot |z-x|^2+ (\xi^*-q^*)\cdot (z_1-x^*)+g(z_1)-g(x^*)\\
&~\leq~\O(1)\cdot \left(1+{1\over 2\bar{\rho}}\right)\cdot |p(0)|\cdot \left(|z-x|^2+|z_1-x^*|^2\right),
\end{split}
\]
and this yields (\ref{P-0}).
\item Otherwise, if $h_z(\tau^*)< \alpha$ then there exists $t_{\beta}\in\, ]\tau^*,{\tau^*_z}]$ such that 
\[
h_z(\tau^*)~<~\alpha~=~V(x)~=~h_z(t_{\beta})\,.
\]
In particular, one has 
\[
\int_{0}^{\tau^*} r\left(y^{x,u^*}(s),u^*(s)\right)ds+g(x^*)~=~\int_{0}^{t_{\beta}} r\left(y^{z,u^*}(s),u^*(s)\right)ds+g(y^{z,u^*}(t_{\beta}))
\]
and this implies  
\[
\begin{split}
\O(1)\cdot |z-x|&~\geq~\int_{\tau^*}^{t_{\beta}} r\left(y^{z,u^*}(s),u^*(s)\right)ds+g(y^{z,u^*}(t_{\beta}))-g(y^{z,u^*}(\tau^*))\\
&~\geq~(t_{\beta}-\tau^*)\cdot r_0-G\cdot |y^{z,u^*}(t_{\beta})-y^{z,u^*}(\tau^*)|~\geq~(t_{\beta}-\tau^*)\cdot(r_0-GN).
\end{split}
\]
Hence, one gets
\bel{t-tau}
t_{\beta}-\tau^*~\leq~{\O(1)\over r_0-GN}\cdot |z-x|\,.
\eeq
Set $z_2\doteq  y^{z,u^*}({{t_\beta}})$. Recalling {(\ref{et-11})}, we estimate 
\[
\begin{split}
- p(0)\cdot (z-x)&~\leq~ \O(1)\cdot |p(0)|\cdot |z-x|^2+ (\xi^*-q^*)\cdot (z_1-x^*)\\ 
&\qquad \qquad\qquad -g(x^*) +h_z(t_{\beta})-\int_{0}^{\tau^*}r(y^{z,u^*}(s),u^*(s))ds\\
&~=~ \O(1)\cdot |p(0)|\cdot |z-x|^2+ (\xi^*-q^*)\cdot (z_2-x^*)+g(z_2)\\ 
&\qquad\qquad\qquad  -g(x^*)+\int_{\tau^*}^{{{t_\beta}}}r(y^{z,u^*}(s),u^*(s))ds+(\xi^*-q^*)\cdot (z_1-z_2)\\
&~\leq~ \O(1)\cdot \left(1+{1\over 2\bar{\rho}}\right)\cdot |p(0)|\cdot \left(|z-x|^2+|z_2-x^*|^2\right)\\
&\qquad\qquad\qquad +\int_{\tau^*}^{{{t_\beta}}}\left[r\big(y^{z,u^*}(s),\omega^*\big)+p^*\cdot f\big(y^{z,u^*}(s),\omega^*\big)\right]ds\\
&~\leq~\O(1)\cdot \left(1+{1\over 2\bar{\rho}}\right)\cdot |p(0)|\cdot \left(|z-x|^2+|{{t_\beta}}-\tau^*|^2+|z_1-x^*|^2\right)\\
&\qquad\qquad\qquad +\int_{\tau^*}^{{{t_\beta}}}\left[r\big(x^*,\omega^*\big)+p^*\cdot f\big(x^*,\omega^*\big)\right]ds.
\end{split}
\]
Thus, (\ref{chos1}) and (\ref{t-tau}) yield (\ref{P-0}).
\end{itemize}
The proof is complete.
\endproof
\medskip

From  Lemma \ref{NN}, Lemma \ref{p1} and the assumption {(\bf H4)}, one immediately derives  an interior sphere property of sublevel set of value function $V$.
\begin{corollary} In addition to  the same assumptions of Lemma \ref{p1},  if {\bf (H4)} holds then the sublevel set $\mathcal{V}(\alpha)$ satisfies an interior sphere condition for every $\alpha>0$.
\end{corollary}

Before proving that $p(0)$ in Lemma \ref{p1} is actually a supergradient of $V$ at $x$, we  provide an exterior sphere condition of $\mathrm{hypo}(V)$ on the boundary of $\mathcal{S}$. 
\begin{lemma}[The regularity of $V$ on $\partial\mathcal{S}$]\label{V-b} Let $x$ be in $\partial\mathcal{S}$. Assume that there exist $q_x\in D^+g(x)$ and $\xi_x\in N^{P}_{\overline{\mathcal{S}^c}}(x)$ realized by a ball of radius $\bar{\rho}>0$ such that $H(x,p_x)=0$ with $p_x=-\xi_{x}+q_x$. Then it holds 
\bel{et1}
-p_x\cdot (z-x)+\beta-g(x)~\leq~ \O(1)\cdot \left(1+{1\over 2\bar{\rho}}\right)\cdot (1+|p_x|)\cdot \left(|z-x|^2+|\beta-g(x)|^2\right)
\eeq
for all $z\in \mathcal{S}^c, \beta\leq V(z)$.
\end{lemma}
{\bf Proof.} Let $u_x\in  U$ be a control such that 
\[
H(x,p_x)~=~-p_xf(x,u_x)-r(x,u_x)~=~0\,.
\]
For any $z\in \mathcal{S}^c$, let $z_1(\cdot)=y^{z,u_x}(\cdot)$ be the trajectory starting from $z$ associated with the constant control $u_x$. For any $\beta\leq V(z)$, two cases are considered:
\v

\n $\bullet$ If $\beta\leq g(z)$ then 
\[
\begin{split}
\beta-g(x)-p_x\cdot (z-x)&~\leq~g(z)-g(x)-p_x\cdot (z-x)\\
&~=~g(z)-g(x)-q_x\cdot (z-x)+ \xi_x\cdot (z-x)\\
&~\leq~\O(1)\cdot \left(1+{1\over 2\bar{\rho}}\right)\cdot (1+|\xi_x|)\cdot |z-x|^2\,.
\end{split}
\]
$\bullet$ Otherwise, assume that  $g(z)<\beta\leq V(z)$. Consider the map $h:[0,\infty[\to\R$ such that 
\[
h(t)~=~\int_{0}^{t}r(z_1(s),u_x)ds+g(z_1(t))\qquad\forall t\geq 0.
\]
Using the condition {\bf (H3)}, we can show as in the proof of Lemma \ref{p1} that the map $t\mapsto h(t)$ is strictly increasing and there exists $t_\beta>0$ such that 
\[
g(z)~=~h(0)~\leq~\beta~=~h(t_{\beta})~\leq~V(z)\,.
\]
Notice that the last inequality implies that $z_1(t_\beta)$ is in the closure of $\mathcal{S}^c$. Otherwise, the first time $t_1$ when $z_1(\cdot)$ touch the target $\mathcal{S}$ is between $[0,t_{\beta}]$ and this yields a contradiction since
\[
h(t_{\beta})~>~h(t_1)~=~\int_{0}^{t_1}r(z_1(s),u_x)ds+g(z_1(t_1))~\geq~V(z)\,.
\]
On the other hand, we can also bound $t_{\beta}$ by 
\bel{b-t-beta}
\begin{split}
t_{\beta}&~\leq~{h(t_{\beta})-h(0)\over { r_0}-GN}~=~{\beta-g(z)\over r_0-GN}~\leq~{1\over r_0-GN}\cdot \left(|\beta-g(x)|+|g(z)-g(x)|\right)\\
&~\leq~{1\over r_0-GN}\cdot \left(|\beta-g(x)|+G\cdot |z-x|\right).
\end{split}
\eeq
Since $\xi_x\in N^{P}_{\overline{\mathcal{S}^c}}(x)$ realized by a ball of radius $\bar{\rho}$, it holds
\[
\begin{split}
\xi_x\cdot (z_1(t_{\beta})-x)&~\leq~{|\xi_x|\over 2\bar{\rho}}\cdot |z_1(t_{\beta})-x|^2\\
&~\leq~\O(1)\cdot \left(1+{1\over 2\bar{\rho}}\right)\cdot |\xi_x|\cdot \left(|\beta-g(x)|^2+|z-x|^2\right).
\end{split}
\]
Thus, we estimate 
\[
\begin{split}
\beta-g(x)&-p_x\cdot (z-x)~=~h(t_{\beta})-g(x)-p_x\cdot (z-z_1(t_\beta)+z_1(t_\beta)-x)\\
&~=~\int_{0}^{t_{\beta}}r(z_1(s),u_x)ds-p_x\cdot (z-z_1(t_\beta))\\
&~~~~~~~~~\qquad+ g(z_1(t_{\beta}))-g(x)-q_x\cdot (z_1(t_\beta)-x)+\xi_x\cdot (z_1(t_\beta)-x) \\
&~\leq~\int_{0}^{t_{\beta}}r(z_1(s),u_x)ds-p_x\cdot (z-z_1(t_\beta))+\O(1)\cdot \left(1+{1\over 2\bar{\rho}}\right)\cdot (1+|\xi_x|)\left|z_1(t_\beta)-x\right|^2.
\end{split}
\]
Notice that 
\[
\left|z_1(t_\beta)-x\right|~\leq~\O(1)\cdot (t_{\beta}+|z-x|)~\leq~\O(1)\cdot \left(|\beta-g(x)|+|z-x|\right)\,,
\]
\[
\begin{split}
-p_x\cdot (z-z_1(t_\beta))&~=~\int_{0}^{t_\beta}p_x\cdot f(z_1(s),u_x)ds~\leq~t_{\beta}\cdot p_x\cdot f(x,u_x)+O(1)\cdot |p_x| \cdot(t^2_{\beta}+t_\beta|z-x|)\\
&~=~t_{\beta}\cdot p_x\cdot f(x,u_x)+O(1)\cdot |p_x| \cdot \left(|\beta-g(x)|^2+|z-x|^2\right),
\end{split}
\]
and 
\[
\begin{split}
\int_{0}^{t_{\beta}}r(z_1(s),u_x)ds&~\leq~t_{\beta}\cdot r(x,u_x)+O(1)\cdot (t^2_{\beta}+t_\beta|z-x|)\\
&~\leq~t_{\beta}\cdot r(x,u_x)+O(1)\cdot  \left(|\beta-g(x)|^2+|z-x|^2\right).
\end{split}
\]
This implies that 
\[
\beta-g(x)-p_x\cdot (z-x)~\leq~\O(1)\cdot \left(1+{1\over 2\bar{\rho}}\right)\cdot (1+|p_x|)\cdot \left(|z-x|^2+|\beta-g(x)|^2\right),
\]
and the proof is complete.
\endproof

\begin{proposition}\label{Pro-1}  Under the same setting in Lemma \ref{p1}, it holds
\[
p(t)~\in~\partial^{P}V(y^{x,u^*}(t))\qquad\forall t\in [0,\tau^*]
\]
and  for all $(z,\beta)\in \mathrm{hypo}(V)$
\bel{hyp-f-1}
{(-p(t),1)\over |(-p(t),1)|}\cdot \big(z-y^{x,u^*}(t),\beta-V(y^{x,u^*}(t))\big)~\leq~\O(1)\cdot \left(1+{1\over 2\bar{\rho}}\right)\cdot \big(|z-y^{x,u^*}(t)|^2+|\beta-V(y^{x,u^*}(t))|^2\big).
\eeq
\end{proposition}
{\bf Proof.} For simplicity,  we prove the statement for $t=0$. Recall that by Lemma \ref{O}, for a.e. $t\in [0,\tau^*]$ it holds
\bel{H=0}
H(p(t),y^{x,u^*}(t))~=~- p(t)\cdot f(y^{x,u^*}(t),u^*(t))-r(y^{x,u^*}(t),u^*(t))~=~0.
\eeq
We set 
\bel{V-x}
\alpha~\doteq~V(x)~=~\int_{0}^{\tau^*}r\left(y^{x,u^*}(s), u^*(s)\right)ds+g(x^*)
\eeq
and, for any given $z\in\mathcal{S}^c$, we denote $\alpha_z\doteq V(z)$ and $\tau_z\doteq\tau^{z,u^*}$. Two main cases are considered:
\medskip

$\bullet$ {\bf Case 1.} Assume that  $\alpha_z\leq\alpha$.  Then we only need to prove (\ref{hyp-f-1}) for $\beta=\alpha_z$, i.e., 
\bel{hyp-f-1-1}
{(-p(0),1)\over |(-p(0),1)|}\cdot (z-x,\alpha_z-\alpha)~\leq~\O(1)\cdot \left(1+{1\over 2\bar{\rho}}\right)\cdot \left(|z-x|^2+|\alpha_z-\alpha|^2\right)
\eeq
There are two subcases:
\begin{itemize}
\item {\bf Subcase 1:} If $g(x^*)\leq \alpha_z$, we have 
\[
\int_{0}^{\tau^*}r(y^{x,u^*}(s),u^*(s))ds~=~\alpha-g(x^*)~\geq~\alpha-\alpha_z.
\]
Thus, by DPP there exists $\tau_1\in [0,\tau^*]$ such that $x_1\doteq y^{x,u^*}(\tau_1)\in \partial\mathcal{V}(\alpha_z)$ and 
\[
V(x_1)~=~\alpha-\int_{0}^{\tau_1}r(y^{x,u^*}(s),u^*(s))ds~=~\alpha_z.
\]
Moreover, from  {\bf (H3)}, it holds 
\bel{est-t1}
\tau_1~\leq~{1\over r_0}\cdot\int_{0}^{\tau_1}r(y^{x,u^*}(s),u^*(s))ds~\leq~{\alpha-\alpha_z\over r_0}.
\eeq
Recalling Lemma \ref{p1}, we obtain
\[
\begin{split}
-p(\tau_1)\cdot (z-x_1)&~\leq~ \O(1)\cdot \left(1+{1\over 2\bar{\rho}}\right)\cdot |p(\tau_1)|\cdot |z-x_1|^2\\
&~\leq~ \O(1)\cdot \left(1+{1\over 2\bar{\rho}}\right)\cdot |p(0)|\cdot \big(|x-z|^2+\tau_1^2\big),
\end{split}
\]
and (\ref{est-t1}) implies
\bel{td-0}
-p(0)\cdot (z-x_1)~\leq~\O(1)\cdot \left(1+{1\over 2\bar{\rho}}\right)\cdot (|p(0)|+1)\cdot \big(|x-z|^2+|\alpha-\alpha_z|^2\big).
\eeq
On the other hand, from (\ref{H=0}), one has  
\[
\begin{split}
-p(0)\cdot (x_1-x)&~=~-\int_{0}^{\tau_1}p(0)\cdot f(y^{x,u^*}(s),u^*(s))ds\\
&~\leq~\O(1)\cdot |\tau_1|^2-\int_{0}^{\tau_1}p(s)\cdot f(y^{x,u^*}(s),u^*(s))ds\\
&~=~\O(1)\cdot |\tau_1|^2+\int_{0}^{\tau_1}r(y^{x,u^*}(s),u^*(s))ds~=~\O(1)\cdot |\tau_1|^2+\alpha-\alpha_z.
\end{split}
\]
This, together with (\ref{est-t1}) and (\ref{td-0}), yields (\ref{hyp-f-1-1}).
\medskip

\item {\bf Subcase 2:}  Otherwise, if $g(x^*)> \alpha_z$ then by DPP we have
\bel{est-t*}
\tau^*~\leq~{1\over r_0}\cdot\int_{0}^{\tau^*}r(y^{x,u^*}(s),u^*(s))ds~=~{\alpha-g(x^*)\over r_0}~\leq~{\alpha-\alpha_z\over r_0}.
\eeq
Recalling (\ref{H=0}), we estimate 
\bel{t1d}
\begin{split}
-p(0)\cdot (z-x)&~=~-p(0)\cdot (z-x^*)-p(0)\cdot (x^*-x)\\
&~\leq~\O(1)\cdot \tau^*\cdot |z-x^*|-p(\tau^*)\cdot (z-x^*)-\int_{0}^{\tau^*}p(0)\cdot f(y^{x,u^*}(s),u^*(s))ds\\
&~\leq~\O(1)\cdot \big(|\tau^*|^2+|z-x|^2\big)-p^*\cdot (z-x^*)-\int_{0}^{\tau^*}p(s)\cdot f(y^{x,u^*}(s),u^*(s))ds\\
&~\leq~\O(1)\cdot (|\tau^*|^2+|z-x|^2)-p^*\cdot (z-x^*)+\int_{0}^{\tau^*}r(y^{x,u^*}(s),u^*(s))ds\\
&~=~\O(1)\cdot (|\tau^*|^2+|z-x|^2)-p^*\cdot (z-x^*)+\alpha-g(x^*)
\end{split}
\eeq
On the other hand, from Lemma \ref{V-b}, one has 
\[
\begin{split}
-p^*\cdot (z-x^*)+\alpha_z-g(x^*)&~\leq~\O(1)\cdot \left(1+{1\over 2\bar{\rho}}\right)\cdot |p^*|\cdot\big(|z-x^*|^2+|\alpha_z-g(x^*)|^2\big)\\
&~\leq~\O(1)\cdot \left(1+{1\over 2\bar{\rho}}\right)\cdot |p(0)|\cdot \left(|z-x|^2+|\tau^*|^2+|\alpha_z-\alpha|^2\right),
\end{split}
\]
and (\ref{t1d}) yields (\ref{hyp-f-1-1}).
\medskip
\end{itemize}

$\bullet$ {\bf Case 2.} If $\alpha_z>\alpha$ then $z\in\mathcal{V}(\alpha)^c$ and by Lemma \ref{p1}, we have 
\[
-p(0)\cdot (z-x)~\leq~\O(1)\cdot |z-x|^2.
\]
In particular, (\ref{hyp-f-1}) holds for all $\beta\leq \alpha$. Thus, we shall prove (\ref{hyp-f-1}) for   $\alpha<\beta\leq \alpha_z$. Three subcases are considered:
\begin{itemize}
\item {\bf Subcase 1}:  If  $\tau_z\leq \tau^*$ then (\ref{tau-ve-12}) holds. Following the argument of case 1 in the proof of Lemma \ref{p1}, set 
\[
z^*= y^{z,u^*}(\tau_z)\in\partial\mathcal{S},\qquad \beta_1~=~\beta-\int_{0}^{\tau_z} r(y^{z,u^*}(s),u^*(s))ds,
\]
we estimate 
\[
\begin{split}
-p(0)\cdot (z-x)+\beta-\alpha&~\leq~-p(\tau^*)\cdot (z^*-x^*)+\O(1)\cdot |z-x|^2+\beta-\alpha\\
&\qquad+\int_{0}^{\tau^*} r(y^{x,u^*}(s),u^*(s))ds-\int_{0}^{\tau_z} r(y^{z,u^*}(s),u^*(s))ds\\
&~=~\O(1)\cdot |z-x|^2-p(\tau^*)\cdot (z^*-x^*)+\beta_1-g(x^*).
\end{split}
\]
Noticing that $\beta_1\leq V(z^*)$ and 
\[
\begin{split}
|\beta_1-g(x^*)|&~\leq~|\beta-\alpha|+\left|\int_{0}^{\tau^*}r(y^{x,u^*}(s),u^*(s))ds-\int_{0}^{\tau_z}r(y^{z,u^*}(s),u^*(s))ds\right|\\
&~\leq~|\beta-\alpha|+ \O(1)\cdot \left(|z-x|+|\tau^*-\tau_z|\right),
\end{split}
\]
we  apply Lemma \ref{V-b} to obtain
\[
\begin{split}
-p(0)\cdot (z-x)+\beta-\alpha&~\leq~\O(1)\cdot \left(|z-x|^2+|z^*-x^*|^2+|\beta_1-g(x^*)|^2\right)\\
&~\leq~\O(1)\cdot \left(|z-x|^2+|\tau^*-\tau_z|^2\right).
\end{split}
\]
Thus, (\ref{tau-ve-12}) yields (\ref{hyp-f-1}) for $t=0$.
\medskip

\item  {\bf Subcase 2}: If  $\tau_z> \tau^*$ and there exists $\tau_1\in [0,\tau_z]$ such that 
\bel{eqc}
\int_{0}^{\tau_1}r(y^{z,u^*}(s),u^*(s))ds~=~\beta-\alpha.
\eeq
Then, we have 
\[
\tau_1~\leq~{\beta-\alpha\over r_0},\qquad z_1~\doteq~y^{z,u^*}(\tau_1)~\in~\overline{\mathcal{V}(\alpha)^c},
\]
and Lemma \ref{p1} yields 
\bel{td-2}
\begin{split}
-p(0)\cdot (z_1-x)&~\leq~\O(1)\cdot \left(1+{1\over 2\bar{\rho}}\right)\cdot |p(0)|\cdot |z_1-x|^2\\
&\leq~\O(1)\cdot \left(1+{1\over 2\bar{\rho}}\right)\cdot |p(0)|\cdot \big(|z-x|^2+|\beta-\alpha|^2\big).
\end{split}
\eeq
From (\ref{H=0}) and (\ref{eqc}), one has 
\[
\begin{split}
-p(0)\cdot (z-z_1)&~=~p(0)\cdot \int_{0}^{\tau_1}f(y^{z,u^*}(s),u^*(s))ds\\
&~\leq~\O(1)\cdot \big(|\tau_1|^2+|z-x|^2\big)+\int_{0}^{\tau_1}p(s)\cdot f(y^{x,u^*}(s),u^*(s))ds\\
&~\leq~\O(1)\cdot \big(|\beta-\alpha|^2+|z-x|^2\big)-\int_{0}^{\tau_1}r(y^{x,u^*}(s),u^*(s))ds\\
&~=~\O(1)\cdot \big(|\beta-\alpha|^2+|z-x|^2\big)-(\beta-\alpha),
\end{split}
\]
and (\ref{td-2}) yields (\ref{hyp-f-1}).
\medskip

\item  {\bf Subcase 3}: Otherwise, if  $\tau_z> \tau^*$ and (\ref{eqc}) does not hold, then by DPP 
\[
\tau^*~\leq~\tau_z~\leq~ {1\over r_0} \cdot \int_{0}^{\tau_z}r(y^{z,u^*}(s),u^*(s))ds~<~{\beta-\alpha\over r_0}.
\]
Set $z^*\doteq y^{z,u^*}(\tau^*)$, we estimate 
\[
\begin{split}
-p(0)\cdot (z-x)&~=~-p(\tau^*)\cdot (z^*-x^*)+\int_{0}^{\tau^*}{d\over ds}\left[ p(s)\cdot \big( y^{z,u^*}(s)-y^{x,u^*}(s)\big)\right]ds\\
&~\leq~\O(1)\cdot |z-x|^2-p(\tau^*)\cdot (z^*-x^*)+\alpha-g(x^*)-\int_{0}^{\tau^*}r(y^{z,u^*}(s),u^*(s))ds,
\end{split}
\]
and Lemma \ref{V-b} implies 
\[
\begin{split}
-p(0)\cdot (z-x)&+\beta-\alpha\\
&~\leq~\O(1)\cdot |z-x|^2-p(\tau^*)\cdot (z^*-x^*)+\beta-\int_{0}^{\tau^*}r(y^{z,u^*}(s),u^*(s))ds-g(x^*)\\
&\leq~\O(1)\cdot \left(|z-x|^2+|z^*-x^*|^2+\left|\beta-\int_{0}^{\tau^*}r(y^{z,u^*}(s),u^*(s))ds-g(x^*)\right|^2\right)\\
&\leq~\O(1)\cdot \left(|z-x|^2+|\beta-\alpha|^2+|\tau^*|^2\right)~\leq~\O(1)\cdot \left(|z-x|^2+|\beta-\alpha|^2\right).
\end{split}
\]
\end{itemize}
The proof is complete.
\endproof
\medskip

As a consequence of Proposition \ref{Pro-0} and Proposition \ref{Pro-1}, we immediately obtain:

\begin{corollary}\label{d-point} Under the same setting in Proposition \ref{Pro-1}, if $V$ is differentiable at $x\in \mathcal{S}^c$ then 
\bel{eq-DV}
\{DV(x)\}~=~\partial^{P}V(x)=N_1^{\dagger}(x).
\eeq
\end{corollary}
{\bf Proof.} Assume that $V$ is differentiable at $x\in \mathcal{S}^c$. Let $u^*$ be an optimal control steering $x$ to $x^*\in\partial\mathcal{S}$.
 By Proposition \ref{Pro-0}, if $N_0(x^*)$ is nonempty then $\partial^{\infty}V(x)$ is nonempty and this contradicts the differentiability of $V$ at $x$. Thus, from Lemma \ref{NN} and {\bf (H4)},  $N_1(x^*)$ is nonempty. Thus, by Proposition \ref{Pro-1}, we have 
\[
\varnothing~\neq~N_1^{\dagger}(x)~\subseteq~\partial^{P}V(x),
\]
and the  differentiability of $V$ at $x$ yields (\ref{eq-DV}).
\endproof

\subsection{Proof of Theorem \ref{Main1} and Theorem \ref{Re-pre}} 
In this subsection, we shall  prove our main theorems by using the previous lemmas and propositions. 
\medskip

{\bf Proof of Theorem \ref{Main1}.} Given $x\in \mathcal{R}\backslash\mathcal{S}$, let  $u^*(\cdot)$ be an optimal control which steers $x$ to $x^*=y^{x,u^*}(\tau^*)\in\partial\mathcal{S}$ with $\tau^*\doteq\tau^{x,u^*}$. By ({\bf H4}) and Lemma \ref{NN}, there exists a proximal normal vector  $\xi^*\in N_0(x^*)\cup N_{1}(x^*)$ to $\overline{\mathcal{S}^c}$ which is realized by a ball of radius  $\rho_0$. 
\begin{itemize}
\item If $\xi^*\in N_0(x^*)$ then let $p(\cdot)$ be the solution of (\ref{adjoint-2}). By Proposition \ref{Pro-0}, the vector $(-p(0),0)\in N^P_{\mathrm{hypo}(V)}(x,V(x))$ is realized by a ball of radius $r(\rho_0)=\O(1)\cdot\ds{ \rho_0\over  1+2\rho_0}>0$.

\item Otherwise,  if $\xi^*\in N_1(x^*)$ then  let $p(\cdot)$ be the solution of (\ref{adjoint-1}) with $p^*=q^*-\alpha \xi^*$ with $q^*\in D^+g(x^*)$, $\alpha>0$ such that $H(x^*,p^*)=0$. By Proposition \ref{Pro-1}, the vector $(-p(0),-1)\in N^P_{\mathrm{hypo}(V)}(x,V(x))$ is realized by a ball of radius $r(\rho_0)=\O(1)\cdot\ds{ \rho_0\over  1+2\rho_0}>0$.
\end{itemize}
In particular, $\mathrm{hypo}(V)$ satisfies a $r(\rho_0)$-exterior sphere condition. Furthermore, from Propositions \ref{Pro-0} and \ref{Pro-1}, it holds
\[
N_0^{\dagger}(x)~\subseteq~\partial^{\infty} V(x),\qquad N_1^{\dagger}(x)~\subseteq~\partial^{P} V(x).
\]
Thus, by the convexity of $\partial^{\infty} V(x)$ and $\partial^{P} V(x)$, and $\partial^{\infty} V(x)+\partial^{P} V(x)= \partial^{P} V(x)$,  one  achieves the inclusions in (\ref{sup-V}).
\endproof
\v

{\bf Proof of Theorem \ref{Re-pre}.} If $x\in\mathrm{int}\mathcal{R}\backslash\mathcal{S}$ is such that $N^P_{\mathrm{hypo}(V)}(x,V(x))$ is pointed, then by Proposition \ref{Pres-f} one has 
\bel{tre}
\partial^{P} V(x)~=~\partial^{\infty} V(x)+\mathrm{co}\big[\partial^{*}V(x)\big],\qquad \partial^{\infty} V(x)~=~\mathrm{co}\big[ \partial^{*,\infty}V(x)\cup \{0\}\big].
\eeq
{\bf 1.} Let us prove (i). Thanks to Theorem \ref{Main1} and (\ref{tre}), in order to conclude it will be sufficient to show that $\partial^{*}V(x)\subseteq N_1^{\dagger}(x)$. Given a vector $p_0\in\partial^{*}V(x)$, there exists a sequence $(x_n)_{n\geq1}$ in $\mathcal{S}^c$ converging to $x$ such that $V$ is differentiable at $x_n$ for every $n\geq1$ and $\ds \lim_{n\to\infty} DV(x_n)=p_0$. By Corollary \ref{d-point}, we have that $\{DV(x_n)\}=N_1^{\dagger}(x_n)$. For every $n\geq1$, let $u^*_n\in\mathcal{U}_{x_n}$ be an optimal control steering $x_n$ to $x^*_n\in\partial\mathcal{S}$ by the trajectory $y_n(\cdot)\doteq y^{x_n,u_n^*}(\cdot)$. Given $q^*_n\in D^+g(x_n^*)$ and $\xi^*_n\in N_1(x_n^*)\bigcap N^P_{\overline{\mathcal{S}^c}}(x^*_n) $ realized by a ball of radius $\rho_0$ and such that $H(x^*_n,q^*_n-\xi^*_n)=0$, let $p_n(\cdot)$ be the unique solution to  
\bel{p-n-t}
p'(t)~=~-D^{T}_xf(y_n(t),u^*_n(t))\cdot p(t)-D_xr(y_n(t),u^*_n(t)),\qquad p(\tau^*_n)~=~q_n^*-\xi_n^*,
\eeq
with $\tau_n^*~\doteq~\tau^{x_n,u_n^*}$.
By Proposition \ref{Pro-1} and Lemma \ref{O}, we have $DV(x_n)=p_n(0)$,
$$(-p_n(t),1)\in N^{P}_{\mathrm{hypo}(V)}(y_n(t),V(y_n(t)))$$ is realized by a ball of radius $\ds\O(1)\cdot {\rho_0\over 1+\rho_0}$ for every $t\in [0,\tau_n^*]$ and 
 \bel{con-1}
-p_n(t)\cdot f(y_n(t),u^*_n(t))- r(y_n(t),u^*_n(t))~=~ H(y_n(t),p_n(t))~=~0\quad a.e.~t\in [0,\tau_n^*].
 \eeq
By Gronwall's inequality, we have 
\[
|p_n(t)|~\leq~\O(1)\cdot |p_n(0)|~=~\O(1)\cdot |DV(x_n)|~\leq~\O(1)\cdot (1+|p_0|)\quad\forall n\geq 1.
\]
Without loss of generality we can assume that $(y_n(\cdot),p_n(\cdot))$ converges uniformly to  $(\bar{y}(\cdot),\bar{p}(\cdot))$ such that  $\bar{y}(\cdot)$ optimally steers $x$ to $x^*\doteq \bar{y}(\tau^*)\in \partial\mathcal{S}$ with $\ds\tau^*=\ds \lim_{n\to\infty}\tau_n^*$, $q^*=\ds \lim_{n\to\infty}q_n^*\in D^+g(x^*)$ and $\xi^*=\ds\lim_{n\to\infty}\xi_n^*\in N^P_{\overline{\mathcal{S}^c}}(x^*)$ is realized by a ball of radius $\rho_0$. For every $t\in\, ]0,\tau^*[ $,  noticing  that  $(-\bar{p}(t),1)=\ds\lim_{n\to\infty} (-p_n(t),1)\in N^{P}_{\mathrm{hypo}(V)}(\bar y(t),V(\bar y(t)))$ and 
\[
V(x(t-\ve))-V(x(t))~\geq~\ve\cdot r_0\quad\forall \ve \in \,]0,\tau^*-t], 
\]
we have  that $\bar{p}(t)\neq 0$. From \cite[Theorem 3.4.4]{PC}, the $C^1$-smoothness of $H$ on $\R^d\times \big(\R^d\backslash\{0\}\big)$ implies  that for every $(z,q)\in \R^d\times \big(\R^d\backslash\{0\}\big)$, it holds
\[
D_pH(z,q)~=~-f(z,u(z,q)),\qquad D_xH(z,q)~=~-D^{T}_xf(z,u(z,q))\cdot q-D_x r(z,u(z,q)),
\]
with $u(z,q)$ being any element of $U$ such that 
\[
H(z,q)~=~-f(z,u(z,q))\cdot q -r(z,u(z,q)).
\]
Thus, $(y_n(\cdot),p_n(\cdot))$ solves the system of ODEs
\bel{DDP-1}
{d\over dt} (x(t),p(t))~=~\big(-D_pH(x(t),p(t)), D_xH(x(t),p(t))\big)
\eeq
By the continuity of $(x,p)\mapsto \big(D_xH(x,p), D_pH(x,p)\big)$ on $\R^d\times \big(\R^d\backslash\{0\}\big)$, we obtain that $(\bar{y}(\cdot),\bar{p}(\cdot))$ is a solution to (\ref{DDP-1}) with 
\[
x(\tau^*)~=~x^*,\qquad p(\tau^*)~=~q^*-\xi^*.
\]
Let $F:[0,\tau^*]\mapsto U $ be such that 
\[
F(t)~=~\{\omega\in U: -f(\bar{y}(t),\omega)\cdot \bar{p}(t) -r(\bar{y}(t),\omega)=H(\bar{y}(t),\bar{p}(t))\}
\]
By the continuity of $r,f$ and $H$, the multifunction $F$ is closed and measurable.  From a standard measurable selection theorem (see e.g. in \cite[Theorem 5.3]{CLSW}), $F$ admits a measurable selection $\bar{u}(t)$. In this case, $\bar{y}(\cdot)=y^{x,\bar{u}}(\cdot)$ and  $\bar{p}(\cdot)$ is the unique solution to the associated adjoint equation
\[
\dot{p}(t)~=~-D^{T}_xf(\bar{y}(t),\bar{u}(t)) \cdot p(t)-D_xr(\bar{y}(t),\bar{u}(t)),\qquad p(\tau^*)~=~q^*-\xi^*.
\]
Since $H(x^*,p^*(\tau^*))=\ds\lim_{n\to\infty}H(x_n^*,p_n^*)=0$, we then obtain that $p_0=\ds\lim_{n\to\infty}p_n(0)=\bar{p}(0)\in N_1^{\dagger}(x)$. Therefore, $\partial^{*}V(x)\subseteq N_1^{\dagger}(x)$ and 
\[
\partial^{P} V(x)~\subseteq~\partial^{\infty} V(x)+\mathrm{co} \big[N_1^{\dagger}(x)\big].
\]
This yields (\ref{sup-V-e}).
\medskip

{\bf 2.} Similarly, assuming that $H^0$ is $C^{1}$ on  $\R^d\times \big(\R^d\backslash\{0\}\big)$, we will show that $\partial^{*,\infty}V(x)\subseteq N_0^{\dagger}(x)$. This implies
\[
 \partial^{\infty} V(x)~=~\mathrm{co}\big[ \partial^{*,\infty}V(x)\cup \{0\}\big]~\subseteq~\mathrm{co}\big[N_0^{\dagger}(x)\cup \{0\}\big]
\]
and (\ref{sup-V-ii}) in the case where $N^P_{\mathrm{hypo}(V)}(x,V(x))$ is pointed. 

Given a unit  vector $ p_0\in\partial^{*,\infty}V(x)\backslash \{0\}$, there exist $x_n\in \mathcal{S}^c$ converging to $x$ such that $V$ is differentiable at $x_n$ and 
\[
p_0~=~\lim_{n\to\infty}{DV(x_n)\over |DV(x_n)|},\qquad \lim_{n\to\infty} |DV(x_n)|~=~+\infty.
\]
Let $u^*_n\in\mathcal{U}_{x_n}$, $x^*_n\in\partial\mathcal{S}$, $\tau^*_n$, $q^*_n\in D^+g(x_n^*)$ and $\xi^*_n\in N_1(x_n^*)$ be as in the previous case and $p^*_n\doteq q^*_n-\xi^*_n$. We suppose  without loss of generality that $y_n(\cdot)\doteq y^{x_n,u^*_n}(\cdot)$ converges uniformly to an optimal trajectory  $\bar{y}(\cdot)$ steering $x$ to $x^*\doteq \bar{y}(\tau^*)\in \partial\mathcal{S}$, with $\ds\tau^*\doteq\ds \lim_{n\to\infty}\tau_n^*$  and $\xi^*\doteq\ds\lim_{n\to\infty}{\xi_n^*\over |\xi_n^*|}\in N^P_{\overline{\mathcal{S}^c}}(x^*)$ realized by a ball of radius $\rho_0$.  Moreover, from (\ref{p-n-t}), one has 
\[
\lim_{n\to\infty}\big|p_n^*\big|~=~\lim_{n\to\infty}\big| p_n(\tau_n^*)|~\geq~\O(1)\cdot \lim_{n\to\infty} |p_n(0)|~=~\O(1)\cdot \lim_{n\to\infty} |DV(x_n)|~=~+\infty,
\]
and  $\tilde{p}_n(\cdot)\doteq{p}_n(\cdot)/|p_n^*|$ solves 
\[
\tilde{p}_n'(t)~=~-D^{T}_xf(y_n(t),u_n(t)) \tilde{p}_n(t)-{D_xr(y_n(t),u_n(t))\over |p_n^*|},\qquad \tilde{p}_n(\tau^*_n)~=~{p_n^*\over |p_n^*|}.
\]
Suppose without loss of generality that $\tilde{p}_n(\cdot)$ converges uniformly to $\tilde{p}(\cdot)$ such that $\tilde{p}(t)\neq 0$ for all $t\in [0,\tau^*]$. We claim that $(\bar{y}(\cdot),\tilde{p}(\cdot))$ solves the system of ODEs
\bel{DDP-0}
{d\over dt} (\bar{y}(t),\tilde{p}(t))~=~\big(-H^0_p(\bar{y}(t),\tilde{p}(t)), H^0_x(\bar{y}(t),\tilde{p}(t))\big).
\eeq
Indeed, from (\ref{con-1}), one has 
\[
\begin{split}
-f(\bar{x}(t),u_n(t))&\cdot \tilde{p}(t)~=~-f(y_n(t),u_n(t))\cdot {\tilde{p}_n(t)}+\O(1)\cdot \big(|y_n(t)-\bar{x}(t)|-|\tilde{p}_n(t)-\tilde{p}(t)|\big)\\
&~\geq~\max_{\omega\in U}\left\{-f(y_n(t),\omega)\cdot {\tilde{p}_n(t)}\right\}+\O(1)\cdot \left(|y_n(t)-\bar{x}(t)|+|p_n(t)-\tilde{p}(t)|+{1\over |p_n^*|}\right)\\
&~=~\max_{\omega\in U}\left\{-f(\bar{y}(t),\omega)\cdot {\tilde{p}(t)}\right\}+\O(1)\cdot \left(|y_n(t)-\bar{x}(t)|+|p_n(t)-\tilde{p}(t)|+{1\over |p_n^*|}\right)\\
&~=~H^0(\bar{y}(t),\tilde{p}(t))+\O(1)\cdot \left(|y_n(t)-\bar{x}(t)|+|p_n(t)-\tilde{p}(t)|+{1\over |p_n^*|}\right),
\end{split}
\]
and this yields 
\bel{ke1}
\lim_{n\to\infty}-f(\bar{y}(t),u_n(t))\cdot \tilde{p}(t)~=~H^0(\bar{y}(t),\tilde{p}(t)).
\eeq
From \cite[Theorem 3.4.4]{PC}, the $C^1$-smoothness of $H^0$ on $\R^d\times \big(\R^d\backslash\{0\}\big)$ implies  that for every $(z,q)\in \R^d\times \big(\R^d\backslash\{0\}\big)$, 
\[
D_pH^0(z,q)~=~-f(z,u(z,q)),\qquad D_xH^0(z,q)~=~-D^{T}_zf(z,u(z,q))\cdot q,
\]
with $u(z,q)$ being any element of $U$ such that 
\[
H^0(z,q)~=~-f(z,u(z,q))\cdot q .
\]
By a contradiction argument, one can derive from (\ref{ke1}) that 
\[
\lim_{n\to\infty}-f(\bar{y}(t),u_n(t))~=~D_pH^0(\bar{y}(t),\tilde{p}(t)),\qquad \lim_{n\to\infty}-D^{T}_xf(\bar{y}(t),u_n(t))\tilde{p}(t)~=~-D_xH^0(\bar{y}(t),\tilde{p}(t)),
\]
and this yields (\ref{DDP-0}). As in the previous step, one can select a measure function $\bar{u}$ such that $\bar{y}(\cdot)=y^{x,\bar{u}}(\cdot)$ and  $\tilde{p}(\cdot)$ is the unique solution to the associated adjoint equation
\[
\tilde{p}'(t)~=~-D^{T}_xf(\bar{y}(t),\bar{u}(t)) \cdot \tilde{p}(t),\qquad \tilde{p}(\tau^*)~=~\xi^*.
\]
On the other hand, recalling (\ref{con-1}), we have 
\[
\begin{split}
-\xi^*\cdot f(x^*,\omega)&~=~\lim_{n\to\infty} {-p_n^*\over |p_n^*|}\cdot f(x^*_n,\omega)\\
&~=~\lim_{n\to\infty}{1\over |p^*_n|}\big[-p_n^*\cdot  f(x^*_n,\omega)-r(x^*_n,\omega)\big]\\
&~\leq~\limsup_{n\to\infty}~{H(p_n^*,x_n^*)\over |p_n^*|}~=~0\qquad\forall \omega\in U.
\end{split}
\]
Since $x^*$ is in the  reachable boundary $\partial\mathcal{S}^{*}$ of target $\mathcal{S}$,  Lemma \ref{NN}  yields  $\xi^*\in N_0(x^*)$. In particular, we have 
\[
p_0~=~\lim_{n\to\infty}{\tilde{p}_n(0)\over |\tilde{p}_n(0)|}~=~{ \tilde{p}(0)\over |\tilde{p}(0)|}~\in~N^{\dagger}_0(x),
\]
and this concludes the proof. 
\endproof
\section{Optimality conditions}\label{sec4}
\setcounter{equation}{0}
In this section, we shall analyze some optimality conditions for the trajectories of our control problem. More precisely, we shall give a complete version of  the Pontryagin's maximum principle by using Proposition \ref{Pro-0} and Proposition \ref{Pro-1}. Then, we will provide sufficient conditions for optimality and investigate the uniqueness of optimal trajectories. 
\begin{proposition}[Maximum principle]\label{Maxi} In addition to {\bf (H1)}-{\bf (H3)},  assume that $g$ is locally semiconcave and $V$ is continuous. Given $x\in\mathcal{R}\backslash\mathcal{S}$, let  $u^*(\cdot)$ be an optimal control which steers $x$ to $x^*\doteq y^{x,u^*}(\tau^*)\in\partial\mathcal{S}$ with $\tau^*\doteq\tau^{x,u^*}$. 
\begin{itemize}
\item [(i).] Let $q^*\in D^+g(x^*)$ and $\xi^*\in N_1(x^*)$ be such that $H(x^*,p^*)=0$ with $p^*=q^*-\xi^*$. Calling $p(\cdot)$ the solution of the ODE
\bel{adj1}
p'(t)~=~-D^{T}_xf(y^{x,u^*}(t),u^*(t))\cdot p(t)-D_xr(y^{x,u^*}(t),u^*(t)),\qquad p(\tau^*)~=~p^*,
\eeq
then for a.e. $t\in [0,\tau^*]$ it holds 
\bel{MP-1}
H(y^{x,u^*}(t), p(t))~=~-p(t)\cdot f\big(y^{x,u^*}(t),u^*(t)\big)-r\big(y^{x,u^*}(t),u^*(t)\big)~=~0.
\eeq
\item [(ii).] Let $\xi^*\in N_0(x^*)$ and $p(\cdot)$ be the solution of the ODE
\bel{adj0}
p'(t)~=~-D^{T}_xf(y^{x,u^*}(t),u^*(t))\cdot p(t),~\qquad p(\tau^*)~=~-\xi^*,
\eeq
then for a.e. $t\in [0,\tau^*]$ it holds 
\bel{MP-0}
H^0(y^{x,u^*}(t), p(t))~=~-p(t)\cdot f\big(y^{x,u^*}(t),u^*(t)\big)~=~0.
\eeq
\end{itemize}
\end{proposition}
{\bf Proof.} {\bf 1.} A proof for (i) is given in Lemma \ref{O}. Here, we shall provide a different proof by using Proposition \ref{Pro-1}. Fix $t\in \,]0,\tau^*[$ a  Lebesgue point of the maps $s\mapsto f(y^{x,u^*}(s),u^*(s))$ and $s\mapsto r(y^{x,u^*}(s),u^*(s))$. From Proposition \ref{Pro-1} $p(t)\in \partial^{P}V(y^{x,u^*}(t))$, so that  for all $(z,\beta)\in\mathrm{hypo}(V)$, it holds 
\bel{Mx-x-1}
-p(t)\cdot \big (z-y^{x,u^*}(t))+ \beta-V\big(y^{x,u^*}(t)\big)~\leq~\O(1)\cdot \left(\big|z-y^{x,u^*}(t)\big|^2+\big|\beta-V\big(y^{x,u^*}(t)\big)\big|^2\right).
\eeq
In particular, choosing $(z,\beta)=\big(y^{x,u^*}(t\pm \ve), V(y^{x,u^*}(t\pm \ve))\big)$ for $\ve>0$  small, we obtain 
\[
\begin{split}
-p(t)\cdot \big [y^{x,u^*}(t\pm\ve)&-y^{x,u^*}(t)\big]+V\big(y^{x,u^*}(t\pm\ve)\big)-V\big(y^{x,u^*}(t)\big)\\
&\leq~\O(1)\cdot \left(\big|y^{x,u^*}(t\pm\ve)-y^{x,u^*}(t)\big|^2+\big|V\big(y^{x,u^*}(t\pm\ve)\big)-V\big(y^{x,u^*}(t)\big)\big|^2\right)\\
&\leq~\O(1)\cdot \left(\ve^2+\Big|\int_{t}^{t\pm \ve}r(y^{x,u^*}(s),u^*(s))ds\Big|^2\right)~\leq~\O(1)\cdot \ve^2,
\end{split}
\]
and this implies 
\[
-p(t)\cdot {y^{x,u^*}(t\pm\ve)-y^{x,u^*}(t)\over \ve}-{1\over \ve}\int_{t}^{t\pm\ve}r\big(y^{x,u^*}(s),u^*(s)\big)ds~\leq~\O(1)\cdot \ve.
\]
Taking $\ve\to 0+$, we  get   the second equality in (\ref{MP-1}). 
\medskip

For every $\omega\in U$, let $y^{\omega}(\cdot)$ be the trajectory starting from $y^{x,u^*}(t-\ve)$ with a constant control $u\equiv\omega$. By the dynamic programming principle, we have 
\[
\begin{split}
V(y^{\omega}(\ve))&~\geq~V\big(y^{x,u^*}(t-\ve)\big)-\int_{0}^{\ve}r(y^{\omega}(s),\omega)ds\\
&~=~V(y^{x,u^*}(t))-\int_{0}^{\ve}\big[r(y^{\omega}(s),\omega)-r(y^{x,u^*}(t-\ve+s),u(t-\ve+s))\big]ds\\
&~\geq~\beta_1~\doteq~V(y^{x,u^*}(t))-\ve\cdot \big[r(y^{x,u^*}(t),\omega)-r(y^{x,u^*}(t),u^*(t))\big]-o(\ve),
\end{split}
\]
Choosing $z=y^{\omega}(\ve)$ and $\beta=\beta_1$ in (\ref{Mx-x-1}), we derive
\[
-p(t)\cdot {y^{\omega}(\ve)-y^{x,u^*}(t)\over \ve}-\big[r(y^{x,u^*}(t),\omega)-r(y^{x,u^*}(t),u^*(t))\big]~\leq~o(\ve).
\]
Taking $\ve\to 0+$, we  obtain 
\[
-p(t)\cdot \big[f(y^{x,u^*}(t),\omega)-f(y^{x,u^*}(t),u^*(t))\big]-\big[r(y^{x,u^*}(t),\omega)-r(y^{x,u^*}(t),u^*(t))\big]~\leq~0,
\]
and this yields the first equality in (\ref{MP-1}). 
\medskip

{\bf 2.} Let us now prove (ii). Fix $t\in \,]0,\tau^*[$ a  Lebesgue point of the map $s\mapsto f(y^{x,u^*}(s),u^*(s))$. From Proposition \ref{Pro-0} $p(t)\in \partial^{\infty}V(y^{x,u^*}(t))$, so that for all $(z,\beta)\in\mathrm{hypo}(V)$, it holds 
\bel{Mx-x}
-p(t)\cdot \big [z-y^{x,u^*}(t)\big]~\leq~\O(1)\cdot \left(\big|z-y^{x,u^*}(t)\big|^2+\big|\beta-V\big(y^{x,u^*}(t)\big)\big|^2\right).
\eeq
For $\ve>0$ sufficiently small, we have 
\[
-p(t)\cdot \big [y^{x,u^*}(t\pm\ve)-y^{x,u^*}(t)\big]~\leq~\O(1)\cdot \left(\big|y^{x,u^*}(t\pm\ve)-y^{x,u^*}(t)\big|^2+\Big|\int_{t}^{t\pm \ve}r(y^{x,u^*}(s),u^*(s))ds\Big|^2\right),
\]
and this yields 
\[
-p(t)\cdot {y^{x,u^*}(t\pm\ve)-y^{x,u^*}(t)\over \ve}~\leq~\O(1)\cdot \ve.
\]
Taking $\ve\geq 0$, we  obtain   the second equality in (\ref{MP-0}).
\medskip

For every $\omega\in U$, let $y^{\omega}(\cdot)$ be the trajectory starting from  $y^{x,u^*}(t)$ with a constant control $u\equiv\omega$. By the dynamic programming principle, we have 
\[
V(y^{\omega}(\ve))~\geq~V\big(y^{x,u^*}(t)\big)-\int_{0}^{\ve}r(y^{\omega}(s),\omega)ds~\geq~V(y^{x,u^*}(t))-\O(1)\cdot\ve.
\]
Choosing $z=y^{\omega}(\ve)$ and $\beta=V(y^{x,u^*}(t))-\O(1)\cdot\ve\leq V(y^{\omega}(\ve))$ in (\ref{Mx-x}), we derive
\[
-p(t)\cdot {y^{\omega}(\ve)-y^{x,u^*}(t)\over \ve}~\leq~\O(1)\cdot \ve.
\]
Taking $\ve\to 0+$, we get 
\[
-p(t)\cdot  \big[f(y^{x,u^*}(t),\omega)- f(y^{x,u^*}(t),u^*(t))\big]~\leq~0.
\]
Hence, the first equality in (\ref{MP-0}) holds.
\endproof
\medskip

We are providing some sufficient optimality conditions for the trajectories starting from $x\in\mathcal{R}\backslash\mathcal{S}$. For the horizontal case in which $\partial^{\infty}V(x)$ is nontrivial, we need an additional assumption on the dynamics:
\begin{itemize}
\item [{\bf (H5)}.] Assume that  for all $(x,p)\in \R^d\times \big(\R^d\backslash\{0\}\big)$, there exists a unique $u_0(x,p)\in U$ such that 
\[
u_0(x,p)~\in~\mathrm{argmax}_{\omega\in  U} [-f(x,\omega)\cdot p],
\]
and the map $(x,p)\mapsto u_0(x,p)$ is locally Lipschitz continuous  in $\R^d\times \big(\R^d\backslash\{0\}\big)$.
 Moreover, for every $(x,p)\in\R^d\times \big(\R^d\backslash\{0\}\big)$ with $H_0(x,p)=0$, there exists a constant $C>0$ such that 
\bel{T-asp}
\left\{\omega\in U: \left|f(x,\omega)\cdot {p\over |p|}\right|~\leq~\delta\right\}~\subseteq~B(u_0(x,p),C\delta)
\eeq
for every $\delta>0$ sufficiently small.
\end{itemize}

\begin{remark} The condition {\bf (H5)} holds for example if  for any $x\in \R^d$ the set $f(x,U)$ is uniformly convex with a $C^2$ boundary and $f(x,\cdot)$ is
invertible with a locally Lipschitz inverse, uniformly with respect to $x$.
\end{remark}

\begin{proposition}[Optimality conditions]\label{Optimality conditions} Assume that {\bf (H1)} holds and $V$ is continuous. Let $x\in\mathcal{S}^c$ be steered to  $x^*\in\partial\mathcal{S}$ by a control   $u(\cdot)$. Set    $x(\cdot)\doteq y^{x,u}(\cdot)$ and $\tau\doteq\tau^{x,u}$. 
\begin{itemize}
\item [(i).] If  there exists a continuous function $p:[0,\tau]\to \R^d$ such that 
\bel{Op-1}
p(t)~\in~\partial^{P}V(x(t))\quad\forall t\in [0,\tau[,
\eeq
and
\bel{Op-2} 
-p(t)\cdot f\big(x(t),u(t)\big)-r\big(x(t),u(t)\big)~=~0\quad\text{a.e. } t\in [0,\tau],
\eeq
then $x(\cdot)$ is an optimal trajectory.
\item [(ii).] Assume  {\bf (H1)}-- {\bf (H3)},  {\bf (H5)} and   suppose that  $V$ is locally H\"older continuous with exponent $\alpha>1/2$. If  there exists a  function $p:[0,\tau]\to \R^d$ such that 
\bel{Op-10}
p(t)~\in~\partial^{*,\infty}V(x(t))\quad\forall t\in [0,\tau[,
\eeq
and
\bel{Op-20} 
-p(t)\cdot f\big(x(t),u(t)\big)=~0\quad\text{a.e. } t\in [0,\tau],
\eeq
then $x(\cdot)$ is an optimal trajectory.
\end{itemize}

%
%
\end{proposition}
{\bf Proof.}  By the dynamic programming principle, let $t\mapsto Z(t)$ be the increasing map in $[0,\tau]$ such that 
\[
 Z(t)~\doteq~ V(x(t))+\ds\int_{0}^{t}r(x(s),u(s))ds,\qquad t\in [0,\tau].
\]
We  will show that $Z(\cdot)$ has zero right-derivative for all $t\in [0,\tau[$. In this case, the 
continuity of  $Z(\cdot)$ yields that the map $t\mapsto Z(t)$ is constant in $[0,\tau]$ (see for instance \cite[p.3]{CLSW}). Thus,  $x(\cdot)$ is an optimal trajectory.

 {\bf 1.} Assuming that (\ref{Op-1})-(\ref{Op-2}) hold for a continuous function $p:[0,\tau]\to \R^d$. We show that $Z(\cdot)$ has zero right-derivative for all $t\in [0,\tau[$.
 For every $t\in [0,\tau[$ and $h\in [0,\tau-t]$ sufficiently small, we have 
%
%
\[
\begin{split}
 \ds\int_{t}^{t+h}r(x(s),u(s))ds&~=~-\int_{t}^{t+h}p(s)\cdot f(x(s),u(s))ds\\
 &~=~-p(t)\cdot \int_{t}^{t+h}f(x(s),u(s))ds+o(h)\\
 &~=~-p(t)\cdot \big[x(t+h)-x(t)\big]+o(h)
 \end{split}
\]
and
 \[
 -p(t)\cdot \big[x(t+h)-x(t)\big]+V(x(t+h))-V(x(t))~\leq~\sigma\cdot |x(t+h)-x(t)|^2~\leq~\sigma\cdot N^2h^2,
 \]
for some $\sigma>0$. Hence, we obtain
\[
Z(t+h)-Z(t)~=~V(x(t+h))-V(x(t))+\ds\int_{t}^{t+h}r(x(s),u(s))ds~\leq~o(h),
\]
and $Z(\cdot)$ has zero right-derivative for all $t\in [0,\tau[$.
\medskip

{\bf 2.} We are now proving (ii). Fix  $0\leq t< t+h\leq \tau$.
Since $p(t)\in  \partial^{*,\infty}V(x(t))$, there exists $\bar{x}_h\in\mathcal{S}^c$ such that  $V$ is differentiable at $\bar{x}_h$ and 
\bel{cl1}
|\bar{x}_h-x(t)|+V(\bar{x}_h)-V(x(t))+ \left|{p(t)\over |p(t)|}-{DV(\bar{x}_h)\over |DV(\bar{x}_h)|}\right|~\leq~h^2,\qquad |DV(x_k)|~\geq~{1\over h}.
\eeq
Let $(u_h(\cdot),x_h(\cdot),p_h(\cdot))$ be a triple of optimal control,  trajectory, and  adjoint arc for $\bar{x}_h$ with $x_h(0)=\bar{x}_h$ and $p_{h}(0)=DV(\bar{x}_h)$. From proposition \ref{Maxi}, for all $s\in [0,h]$, it holds
\[
\begin{split}
0~=~H(p_h(s),x_h(s))&~=~\max_{\omega\in U}\big\{-p_h(s)\cdot f(x_h(s),\omega)-r(x_h(x),\omega)\big\}\\
&~=~-p_h(s)\cdot f(x_h(s),u_h(s))-r(x_h(x),u_h(s)),
\end{split}
\]
(\ref{cl1}) implies 
\[
\begin{split}
\left|{{p(t)\over |p(t)|}}\cdot f(x(t),u_h(s))\right|&~\leq~\left|{{DV(\bar{x}_h)\over |DV(\bar{x}_h)|}}\cdot f(x_h(s),u_h(s))\right|+\O(1)h\\[3mm]
&~\leq~\left|{{p_h(s)\over |p_h(s)|}}\cdot f(x_h(s),u_h(s))\right|+\O(1)h~=~ \O(1)h.
\end{split}
\]
On the other hand, one also has 
\[
\begin{split}
H_0(x(t),p(t)/|p(t)|)&~=~\max_{\omega\in U}\left\{ - {p(t)\over |p(t)|}\cdot f(x(t),\omega)\right\}\\
&~=~\lim_{h\to 0}\left(\max_{\omega\in U}\left\{ - {DV(\bar{x}_h)\over |DV(\bar{x}_h)|}\cdot f(x(t),\omega)-{r(\bar{x}_h,\omega)\over |DV(\bar{x}_h)|}\right\}+\O(1)\cdot h\right)\\
&~=~\lim_{h\to 0} \left({H(\bar{x}_h,DV(\bar{x}_h))\over |DV(\bar{x}_h)|}+\O(1)\cdot h\right)~=~0.
\end{split}
\]
Thus, from (\ref{T-asp}) and (\ref{Op-20}), one gets
\[
\big|u_h(s)-u(t)\big|~=~\big|u_h(s)-u_0(p(t)/|p(t)|,x(t))\big|~\leq~ \O(1)h,\qquad s\in [0,h].
\]
In particular, 
\[
\begin{split}
\big|\dot{x}_h(s)-\dot{x}(t+s)\big|&~=~\big|f(x_h(s),u_h(s))-f(x(t+s),u(t+s))\big|\\
&~\leq~\O(1)\cdot\left(|x_h(s)-x(t+s)|+|u_h(s)-u(t)|+|u(t)-u(t+s)|\right)~\leq~\O(1) h,
\end{split}
\]
and this yields 
\bel{e-t1-k}
|x_h(s)-x(t+s)|~\leq~\O(1) h^2\qquad\forall s\in [0, h].
\eeq
Since $V$ is locally H\"older continuous with exponent $\alpha>1/2$, it holds
\bel{V-V}
\begin{split}
V(x(t+h))-V(x(t))&~=~V(x(t+h))-V(x_h(h))+V(x_h(h))-V(\bar{x}_h)\\
&~~~~~~~~+V(\bar{x}_h)-V(x(t))~=~-\int_{0}^{h}r(x_h(s),u_h(s))ds+o(h)\\
&~=~-\int_{t}^{t+h}r(x(s),u(s))ds+o(h),
\end{split}
\eeq
and this implies that $Z(\cdot)$ has zero right-derivative for all $t\in [0,\tau[$.
\endproof

\begin{example} 
In Proposition \ref{Optimality conditions}, the conditions \eqref{Op-10} and \eqref{Op-20} alone are not sufficient for the optimality of the trajectory $x(\cdot)$.
Indeed, consider the minimum time problem associated to the set 
$$\mathcal{S}~=~\big\{(x_1,x_2)\in\R^2 : x_2\leq x_1^{5/3}\big\},$$
the function $f((x_1,x_2),u)=(u,0)$ and the control set $U=[-1,1]$. The point $(-1,0)$ is steered to $(0,0)\in\partial\mathcal{S}$ by the constant control $u\equiv1/2$. The vector $p=(0,1)$ satisfies 
$$
p\in\partial^{*,\infty}T_{\mathcal{S}}(y^{x,u}(t))\quad\text{  and }\quad p\cdot f(y^{x,u}(t),u(t))~=~0\qquad\forall t\in [0,2].
$$
Yet, the trajectory $y^{x,u}(\cdot)$ is not optimal. Remark that condition {\bf (H5)} is not safisfied in this example.
\end{example}

To complete this section, we shall  establish necessary and sufficient conditions for the uniqueness of  an optimal control for a given point $x\in \mathcal{R}\backslash\mathcal{S}$. 
%
%

\begin{proposition}\label{uniq} Under the same hypotheses of Theorem \ref{Main1}, assume that both $H$ and $H^0$ are  $C^{1,1}$ on $\R^d\times \big(\R^d\backslash\{0\}\big)$. For every  $x\in \mathcal{R}\backslash\mathcal{S}$, the followings hold:
\begin{itemize}
 \item If $\mathrm{hypo}(V)$   is differentiable at $(x,V(x))$, then there is a unique optimal trajectory steers $x$ to $\mathcal{S}$. 
 \item Moreover, the reverse implication holds provided that the boundary of $\mathcal{S}$ is $C^{1,1}$ and $g$ is differentiable on $\partial\mathcal{S}$.
 \end{itemize}

%
\end{proposition}

{\bf Proof.} {\bf 1.} Assume that $\mathrm{hypo}(V)$  is differentiable at $(x,V(x))$. From Lemma \ref{Diff},  there exists a unit vector $(v,\eta)\in N^{P}_{\mathrm{hypo}(V)}(x,V(x))$ such that
\[
N^{P}_{\mathrm{hypo}(V)}(x,V(x))~=~ \{\lambda\cdot(v,\eta): \lambda\in [0,\infty[\}.
\]
Two cases are considered:
\begin{itemize}
\item If $\eta>0$, then $\partial^{P}V(x)=\Big\{-\dfrac{v}{\eta}\Big\}$ and $\partial^{\infty}V(x)=\{0\}$. By Theorem \ref{Main1}, we have $N^{\dagger}_1(x)=\Big\{\dfrac{v}{\eta}\Big\}$ and $N^{\dagger}_0(x)=\big\{0\big\}$. Let $u^*(\cdot)$ be an optimal control which steers $x$ to $x^*\doteq y^{x,u^*}(\tau^*)\in\partial\mathcal{S}$ with $\tau^*\doteq\tau^{x,u^*}$ and let $p(\cdot)$ be the solution of (\ref{adj1}), with $p^*=q^*-\xi^*$, $q^*\in D^+g(x^*)$ and $\xi^*\in N_1(x^*)$ such that $H(x^*,p^*)=0$. In this case,  $p(t)\neq 0$ for all $\in [0,\tau^*]$. Moreover, by Proposition \ref{Maxi} (i), we have  
\[
H(y^{x,u^*}(t),p(t))~=~ - f(y^{x,u^*}(t), u^*(t))-r(y^{x,u^*}(t),u^*(t))\quad\forall t\in [0,\tau^*].
\]
From \cite[Theorem 7.3.6]{PC}, the $C^{1,1}$-smoothness of  $H$  on $\R^d\times \big(\R^d\backslash\{0\}\big)$ implies  that for every $(z,q)\in \R^d\times \big(\R^d\backslash\{0\}\big)$
\[
D_pH(z,q)~=~-f(z,u(z,q)),\qquad D_xH(z,q)~=~-D^{T}_xf(z,u(z,q))\cdot q-D_xr(z,u(z,q)),
\]
with $u(z,q)$ being any element of $U$ such that 
\[
H(z,q)~=~-f(z,u(z,q))\cdot q -r(z,u(z,q)).
\]
Thus,  $\big(y^{x,u^*}(\cdot),p(\cdot)\big)$ is the unique solution of the Cauchy problem
\[
\begin{cases}
 y'(t)~=~-D_pH\big(y(t),p(t)\big),\\[1mm]
 p'(t)~=~D_xH\big(y(t),p(t)\big),
 \end{cases}\qquad
 \begin{cases}
 y(0)~=~x,\\[1mm]
 p(0)~=~\dfrac{v}{\eta}.
 \end{cases}
 \]

\item  Otherwise, if $\eta=0$ then $\partial^{P}V(x)=\varnothing$ and $\partial^{\infty}V(x)=\{-sv:s\in [0,\infty[\}$.  By Theorem \ref{Main1}, we have $N^{\dagger}_1(x)=\varnothing$ and $N^{\dagger}_0(x)=\{sv:s\in [0,\infty[\}$. Let $u^*(\cdot)$ be an optimal control which steers $x$ to $x^*\doteq y^{x,u^*}(\tau^*)\in\partial\mathcal{S}$ with $\tau^*\doteq\tau^{x,u^*}$ and let  $p(\cdot)$ be the solution of (\ref{adj0}) with $\xi^*\in N_0(x^*)$ such that $\ds\max_{\omega\in U}\{\xi^*\cdot f(x^*,\omega)\}=0$. In this case,  $p(t)\neq 0$ for all $\in [0,\tau^*]$. Moreover, by Proposition \ref{Maxi} (ii), we have  
\[
H^0(y^{x,u^*}(t),p(t))~=~ - f(y^{x,u^*}(t), u^*(t))\cdot p(t)\quad\forall t\in [0,\tau^*].
\]
Again, from \cite[Theorem 7.3.6]{PC}, the $C^{1,1}$-smoothness of  $H^0$ is of class $C^{1,1}$ on $\R^d\times \big(\R^d\backslash\{0\}\big)$ implies  that for every $(z,q)\in \R^d\times \big(\R^d\backslash\{0\}\big)$
\[
D_pH^0(z,q)~=~-f(z,u(z,q)),\qquad D_xH^0(z,q)~=~-D^{T}_xf(z,u(z,q))\cdot q,
\]
with $u(z,q)$ being any element of $U$ such that 
\[
H^0(z,q)~=~-f(z,u(z,q))\cdot q.
\]
Thus,  $\big(y^{x,u^*}(\cdot),p(\cdot)\big)$ is the unique solution of the Cauchy problem
\[
\begin{cases}
 y'(t)~=~-D_pH^0\big(y(t),p(t)\big),\\[1mm]
 p'(t)~=~D_xH^0\big(y(t),p(t)\big),
 \end{cases}\qquad
 \begin{cases}
 y(0)~=~x,\\[1mm]
 p(0)~=~\dfrac{v}{\eta}.
 \end{cases}
 \]
\end{itemize}
In both cases,  we show that  $y^{x,u^*}(\cdot)$ is   a unique optimal trajectory for $x$.
\medskip

{\bf 2.} Assume that there is a unique optimal control $u^*$ that steers $x$ to $x^*\in\partial\mathcal{S}$. Since the boundary of $\mathcal{S}$ is $C^{1,1}$, then $N^P_{\overline{\mathcal{S}^c}}(x^*)=[0,\infty[\cdot\xi^*$, where $\xi^*$ is the standard unit outer normal to $\overline{\mathcal{S}^c}$ at $x^*$.
Two cases are considered:
\begin{itemize}
\item If $\ds\max_{\omega\in U}\{\xi^*\cdot f(x^*,\omega)\}=0$, then $N^P_{\overline{\mathcal{S}^c}}(x^*)=N_0(x^*)$. By Theorem \ref{Re-pre}, $\partial^{*}V(x)\subseteq N_1^{\dagger}(x)=\varnothing$ and $\partial^{*,\infty}V(x)\subseteq N_0^{\dagger}(x)=\,]0,\infty[\,\cdot\, p(0)$, where $p(\cdot)$ is the solution of (\ref{adjoint-22}). By Corollary \ref{1-d-normal}, $\mathrm{hypo}(V)$ is differentiable at $(x,V(x))$.

\item If $\ds\max_{\omega\in U}\{\xi^*\cdot f(x^*,\omega)\}>0$, then $N^P_{\overline{\mathcal{S}^c}}(x^*)=N_1(x^*)$. By Theorem \ref{Re-pre}, $\partial^{*,\infty}V(x)\subseteq N_0^{\dagger}(x)=\varnothing$.  Moreover, calling $\lambda$ the unique positive real number such that $H(x^*,Dg(x^*)-\lambda\xi^*)=0$, we have $\partial^{*}V(x)\subseteq N_1^{\dagger}(x)=\{p(0)\}$,
where $p(\cdot)$ is the solution of (\ref{adjoint-12}), with $p^*=Dg(x^*)-\lambda\xi^*$. Again by Corollary \ref{1-d-normal}, $\mathrm{hypo}(V)$ is differentiable at $(x,V(x))$.
\end{itemize}
The proof is complete.
\endproof
\begin{example}
In Proposition \ref{uniq}, the differentiability of $\mathrm{hypo}(V)$ at $(x,V(x))$ is not sufficient for the uniqueness of the optimal trajectory steering $x$ to $\mathcal{S}$. Indeed, consider the minimum time problem associated to the set 
$$\mathcal{S}~=~\{(x_1,x_2)\in\R^2 : x_2\leq (x_1-1)^3(x_1+1)^3\},$$
the function $f((x_1,x_2),u)=(u,0)$ and the control set $U=[-1,1]$. In this case, the minimum time function is 
\[
T_{\mathcal{S}}(x_1,x_2)~=~ \left(1+x_2^{1/3}\right)^{1/2}-|x_1|,\qquad (x_1,x_2)\in\mathcal{S}^c.
\]
By a direct computation, we have
\[
N^{P}_{\mathrm{hypo}(T_{\mathcal{S}})}(0,0,1)=\{0\}\times\R_-\times\{0\},
\]
and this implies that, $\mathrm{hypo}(T_{\mathcal{S}})$ is differentiable at the point $(0,0)$. However  $(0,0)$ admits two distinct optimal trajectories
\[
x_1(t)~=~(t,0),\qquad  x_2(t)~=~(-t,0),\qquad t\in [0,1].
\]
Remark that $H^0$ is not $C^{1,1}$ on $\R^d\times \big(\R^d\backslash\{0\}\big)$ in this example.
%
%
\end{example}

\begin{remark} The Hamiltonian $H$ is $C^{1,1}$  on $\R^d\times \big(\R^d\backslash\{0\}\big)$ if  for every $(x,p)\in \R^d\times \big(\R^d\backslash\{0\}\big)$ there exists a unique $u^*\in U$ such that
\[
-f(x,u^*)\cdot p~=~\max_{\omega\in U}~\{-f(x,\omega)\cdot p-r(x,\omega)\},
\]
and the map $(x,p)\mapsto u^*(x,p)$ is locally Lipschitz continuous in $\R^d\times (\R^d\backslash\{0\})$.

Similarly, the horizontal Hamiltonian $H^0$ is $C^{1,1}$  on $\R^d\times \big(\R^d\backslash\{0\}\big)$ if  for every $(x,p)\in \R^d\times \big(\R^d\backslash\{0\}\big)$ there exists a unique $u^*\in U$ such that
\[
-f(x,u^*)\cdot p~=~\max_{\omega\in U}~\{-f(x,\omega)\cdot p\},
\]
and the map $(x,p)\mapsto u^*(x,p)$ is locally Lipschitz continuous in $\R^d\times (\R^d\backslash\{0\})$.
\end{remark}
\section{Structure of singular and regular sets}\label{sec5}
\setcounter{equation}{0}
\subsection{$\mathcal{H}^{d-1}$-rectifiable singular set of $V$}  
In this subsection, we show that the singular set $\Sigma_V$ of $V$, defined by 
\bel{Sing-V}
\Sigma_V~\doteq~\left\{x\in\mathcal{S}^c:V~\text{is~not~differentiable~at}~x\right\},
\eeq
is countably $\mathcal{H}^{d-1}$-rectifiable under a suitable assumption on the target.
%
As a consequence, $V$ is  a function of special bounded variation (SBV). In order to do so, let us prove the following lemma related to  the set of points where the proximal horizontal superdifferential of $V$  contains a nonzero vector
\bel{V-infty}
\Sigma_{V,\infty}~=~\left\{x\in \mathcal{S}^c: \partial^{\infty}V(x)\neq \{0\}\right\}.
\eeq
\begin{lemma}\label{Non-lips}
In addition to {\bf (H0)}-{\bf (H4)}, assume that $g$ is locally semiconcave and $V$ is continuous. Then,  $x\in \mathcal{S}^c$ is a  non-Lipschitz point of $V$ if and only if there exists an optimal trajectory steering $x$ to $x^*\in\partial\mathcal{S}$ such that $N_0(x^*)\neq\varnothing$.
\end{lemma}
{\bf Proof.} The implication $(\Leftarrow)$ follows from Proposition \ref{Pro-0}. In order to prove the converse, assume that $V$ is non-Lipschitz at $x\in \mathcal{S}^c$.  By Theorem \ref{Main1}, $\mathrm{hypo}(V)$ satisfies a $\rho(\cdot)$-exterior sphere condition. Thus, Proposition \ref{hypo-upp} implies that there exists a unit proximal horizontal supergradient $\xi\in \partial^{\infty}V(x) $ such that 
\bel{V-infty}
-\xi\cdot (z-x)~\leq~\O(1)\cdot \left(|z-x|^2+|\beta-V(x)|^2\right)\qquad\forall (z,\beta)\in\mathrm{hypo}(V).
\eeq
Since $V$ is differentiable almost everywhere, there exists  a sequence $(x_n)_{n\geq1}$ in $ \mathcal{S}^c$ such that $V$ is differentiable at $x_n$ and 
\bel{x-x1}
\left|x-x_n-{\xi\over n}\right|~\leq~{1\over n^2}\qquad\forall n\geq 1.
\eeq
Choosing $z=x_n$ in (\ref{V-infty}), we obtain
\[
\begin{split}
{1\over n}&~\leq~-\xi\cdot (x_n-x)+{1\over n^2}~\leq~\O(1)\cdot \left(|x_n-x|^2+|\beta-V(x)|^2\right)+{1\over n^2}\\
&~\leq~\O(1)\cdot \left({1\over n^2}+|\beta-V(x)|^2\right)\qquad\forall \beta\leq V(x_n),
\end{split}
\]
and this implies
\bel{v}
V(x)-V(x_n)~=~|V(x_n)-V(x)|~\geq~ {\O(1)\over  n^{1/2}}\qquad\mathrm{for}~n\geq 1~\mathrm{large~enough}.
\eeq
Since $(-DV(x_n),1)\in N^{P}_{\mathrm{hypo}(V)}(x_n,V(x_n))$ is realized by a ball of radius $\rho(x_n)$, one has 
\[
{(-DV(x_n),1)\over |(-DV(x_n),1)|}\cdot (x-x_n, V(x)-V(x_n))~\leq~{1\over \rho(x_n)}\cdot \left(|x-x_n|^2+|V(x)-V(x_n)|^2\right),
\]
and (\ref{x-x1})-(\ref{v}) yield $\ds\lim_{n\to\infty} |DV(x_n)|=+\infty$. As in the proof of Theorem  \ref{Re-pre}, there exist $u^*_n\in\mathcal{U}_{x_n}$, $q^*_n\in D^+g(x_n^*)$, and $\xi^*_n\in N_1(x_n^*)\bigcap N^P_{\overline{\mathcal{S}^c}}(x^*_n) $ realized by a ball of radius $\rho_0$ with 
$$
x_n^*~\doteq~ y^{x_n,u_n^*}(\tau_n^*),\qquad \tau_n^*~\doteq~\tau^{x_n,u_n^*},
$$
 such that 
 \bel{H-nn}
 H(x_n^*,q_n^*-\xi^*_n)~=~0,\qquad \lim_{n\to\infty}|\xi_n^*|~=~+\infty.
 \eeq
 By {\bf (H0)}, we can suppose without loss of generality that the sequence of optimal trajectories $y^{x_n,u^*_n}(\cdot)$ converges uniformly to an optimal trajectory $y^{x,u^*}(\cdot)$ and $\ds\lim_{n\to\infty}{\xi_n^*\over |\xi_n^*|}=\xi^*\in N^P_{\overline{\mathcal{S}^c}}(x^*)$ with $x^*\doteq y^{x,u^*}(\tau^{x,u^*})$. From (\ref{H-nn}), for any $\omega\in U$ one has  
 \[
 \begin{split}
 \xi^*\cdot f(x^*,\omega)&~=~
 \lim_{n\to\infty} {\xi^*_n-q_n^*\over |\xi^*_n-q^*_n|}\cdot f(x_n^*,\omega)~=~\lim_{n\to\infty} \left[-{q_n^*-\xi^*_n\over |\xi^*_n-q^*_n|}\cdot f(x_n^*,\omega)-{{r(x_n^*,\omega)}\over  |\xi^*_n-q^*_n|}\right]\\
 &~\leq~\limsup_{n\to\infty}{H(x_n^*,q_n^*-\xi^*_n)\over  |\xi^*_n-q^*_n|}~=~0,
 \end{split}
 \]
 and this yields 
 \[
 \max_{\omega\in U}~\{\xi^*\cdot f(x^*,\omega)\}~\leq~0.
 \]
 Finally, noticing that  $x^*$ is in the reachable boundary $\partial\mathcal{S}^*$, one can conclude the proof by Lemma \ref{NN}.
%
%
%
%
\endproof
Before stating our main result of the subsection, let us introduce the set
\[
\partial\mathcal{S}^*_{1}~\doteq~\left\{x\in \partial\mathcal{S}^*: \mathrm{dim} \left[N^{P}_{\overline{\mathcal{S}^c}}(x)\right]=1 \right\},
\]
where $\partial\mathcal{S}^*$ is the reachable boundary of $\mathcal{S}$ defined in (\ref{S*}).
\begin{theorem}\label{H1-re} In addition to {\bf (H0)}-{\bf (H4)}, assume that $g$ is locally semiconcave,  $V$ is continuous and $H^0$ is $C^{1,1}$ on $\R^d\times(\R^d\backslash\{0\})$.  For every $x\in \partial\mathcal{S}^*_{1}$, let ${\bf n}(x)$ be the unit proximal normal vector to $\overline{\mathcal{S}^c}$ at $x$. If the set 
\[
\partial\mathcal{S}^{*,0}_{1}~\doteq~\left\{x\in \partial\mathcal{S}^*_1:\max_{\omega\in U}~\{{\bf n}(x)\cdot f(x,\omega)\}=0 \right\}
\]
is countably  $\mathcal{H}^{d-2}$-rectifiable and ${\bf n} $ is locally Lipschitz on $\partial\mathcal{S}^{*,0}_{1}$, then $\Sigma_V$ is countably $\mathcal{H}^{d-1}$-rectifiable.
\end{theorem}
{\bf Proof.} 
{\bf 1.} From Theorem \ref{Main1} and Proposition \ref{hypo-upp}, $V$ is locally semiconcave in the open set $\R^d\backslash (\mathcal{S}\bigcup \Sigma_{V,\infty})$ and this implies that  the set $\Sigma_V\backslash \Sigma_{V,\infty}$ is countably $\mathcal{H}^{d-1}$-rectifiable. Moreover, let us decompose $\Sigma_{V,\infty}=\Sigma^1_{V,\infty}\cup \Sigma^2_{V,\infty}$ with 
\[
\Sigma^1_{V,\infty}~\doteq~\left\{x\in\Sigma_{V,\infty}:\mathrm{dim} \left[N^{P}_{\mathrm{hypo}(V)}(x,V(x))\right]=1 \right\},\]
\[\Sigma^2_{V,\infty}~\doteq~\left\{x\in\Sigma_{V,\infty}:\mathrm{dim} \left[N^{P}_{\mathrm{hypo}(V)}(x,V(x))\right]\geq 2 \right\}.
\]
By \cite[Theorem 1.1]{ANV}, $\Sigma^2_{V,\infty}$ is also countably $\mathcal{H}^{d-1}$-rectifiable. Hence, it remains to  prove  the countable $\mathcal{H}^{d-1}$-rectifiability for $\Sigma^1_{V,\infty}$. By Lemma \ref{Non-lips}, for every given  $x\in \Sigma^1_{V,\infty}$ there exists  a nonzero vector $p(0)\in N^{\dagger}_0(x)$ such that $p(\cdot)$ is the unique solution of the ODE
\bel{adjoint-2}
p'~=~-D^{T}_xf(y^{x,u^*},u^*)\cdot p,\qquad p(\tau^*)~=~-\xi^*,
\eeq
for some  $u^*\in\mathcal{U}_x$ and unit vector $\xi^*\in N_0\big(x^*\big)$, with $x^*\doteq y^{x,u^*}(\tau^*)$ and $\tau^*\doteq \tau^{x,u^*}$. Recalling Proposition \ref{Maxi}, we have 
\bel{H00}
\max_{\omega\in U}~\{\xi^*\cdot f(x^*,\omega)\}~=~0,
\eeq
and
\[
u^*(t)~\in~\argmax_{\omega\in U}~\{- p(t)\cdot f(y^{x,u^*}(t),\omega)\}\qquad a.e.~t\in [0,\tau^*].
\] 
Since $\mathrm{dim} \left[N^{P}_{\mathrm{hypo}(V)}(x,V(x))\right]=1$,  from Theorem  \ref{Main1} one has that  $\mathrm{dim} \left[N^{P}_{\overline{\mathcal{S}^c}}(x^*)\right]=1$. This particularly yields 
\[
x^*~\in~\partial\mathcal{S}^{*,0}_{1},\qquad \xi^*~=~{\bf n}(x^*).
\]
Since the map $(x,p)\mapsto H^0(x,p)\doteq \ds\max_{\omega\in U}\left\{-p\cdot f(x,\omega)\right\} $ is  of class $\mathcal{C}^{1,1}$ in $\R^d\times (\R^d\backslash\{0\})$, one has that 
\[
D_pH^0(x,p)~=~-f(x,u^*(x,p)),\qquad D_xH^0(x,p)~=~-D^{T}_xf(x,u^*(x,p))\cdot p,
\]
with  $u^*(x,p)\in U$ being unique such that
\[
-f(x,u^*)\cdot p~=~\max_{\omega\in U}~\{-f(x,\omega)\cdot p\}.
\]
Moreover, $ \big(y^{x,u^*}(\cdot),p(\cdot)\big)$ solves  the characteristics of the PDE ($H^0(x,\nabla V)=0)$
\bel{Ba-w}
\begin{cases}
\dot{X}(t)~=~-D_pH^0(X(t),P(t)),\\[3mm]
\dot{P}(t)~=~D_xH^0(X(t),P(t)),
\end{cases}
\qquad\qquad 
\begin{cases}
X(\tau^*)~=~x^*,\\[3mm]
P(\tau^*)~=~-{\bf n}(x^*).
\end{cases}
\eeq

{\bf 2.} Consider the map $\phi:[0,\infty[\times \partial\mathcal{S}^{*,0}_{1}\to \R^d$ such that 
\[
\phi(t,z)~=~X(0;t,z)\qquad\forall (t,z)\in [0,\infty[\times \partial\mathcal{S}^{*,0}_{1},
\]
where $(X,P)(\cdot;t,z)$ is the solution to the system of ODEs in (\ref{Ba-w}) with $(X(t),P(t))=(z,-{\bf n}(z))$. By the  locally Lipschitz continuity of ${\bf n}(\cdot)$ in  $\partial\mathcal{S}^{*,0}_{1}$, $\phi$ is  locally Lipschitz in $[0,\infty[\times \partial\mathcal{S}^{*,0}_{1}$. Therefore, the countable  $\mathcal{H}^{d-2}$-rectifiability of $\partial\mathcal{S}^{*,0}_{1}$ implies that the set $\Sigma^1_{V,\infty}\subseteq\phi\big([0,+\infty[\times  \partial\mathcal{S}^{*,0}_{1} \big)$ is countably $\mathcal{H}^{d-1}$-rectifiable. The proof is complete.
\endproof
\medskip

By Theorem \ref{Main1} and Proposition \ref{hypo-upp}, the hypograph of $V$ satisfies an exterior sphere condition and this yields   $V\in BV_{loc}( \mathcal{S}^c)$. In general, functions with hypograph satisfying an exterior sphere condition do not belong to the class of functions of special bounded variation (see \cite[Proposition 4]{GAW1}).  By Theorem \ref{H1-re}, the singular part of $DV(\cdot)$ is concentrated on a set of $\sigma$-finite $\mathcal{H}^{d-1}$-measure. Thus, one obtains the following corollary.

\begin{corollary} Under the same assumptions of Theorem \ref{H1-re},  $V$  is a function of special bounded variation in $ \mathcal{S}^c$.
\end{corollary}

\begin{remark} In the case where the target $\mathcal{S}$ is defined by 
\[
\mathcal{S}~=~\left\{x\in \R^{d}: h(x)~\leq~0\right\}
\]
for some $h\in \mathcal{C}^{2}_b(\R^d)$ with $\nabla h(x)\neq 0$ on $\partial \mathcal{S}$, the set $\partial\mathcal{S}^{*,0}_{1}$ in Theorem \ref{H1-re} can be expressed by 
\[
\partial\mathcal{S}^{*,0}_{1}~=~\big\{x\in \partial\mathcal{S}^*:H^0(x,-\nabla h(x))=0 \big\}.
\]
Assume that for every $x\in  \partial\mathcal{S}^*$ and $\lambda \in\R$, it holds
\[
D_xH^0(x,-\nabla h(x))-D_pH^0(x,-\nabla h(x))\cdot \nabla^2h(x)~\neq~\lambda \nabla h(x).
\]
Then the set $\partial\mathcal{S}^{*,0}_{1}$ is  $\mathcal{H}^{d-2}$-rectifiable. 
\end{remark}
\subsection{Propagation of singularities and differentiablity}
In this subsection, we shall  establish  results about the propagation of singularities and also the differentiability of $V$ along optimal trajectories  under the following additional hypothesis on the Hamiltonian $H$:
\begin{itemize}
\item[ {\bf (H6).}] For any fixed $x\in \R^d$, if $H(x,p)=0$ for all $p$  in a convex set $C$, then $C$ is a
singleton.
\end{itemize}
In this connection, let us recall the definition of  optimal points and provide a necessary condition on  supergradient and horizontal gradient at such points under a pointedness assumption.

\begin{definition}[Optimal points] A point $x\in \mathcal{S}^c$ is called an optimal point of (\ref{Cont})-(\ref{V1}) if there exists an optimal trajectory crossing $x$ at time $t>0$.
\end{definition}
\begin{lemma}\label{H-optimal} Under the same assumptions of Theorem \ref{Main1}, for every optimal point $x\in\mathcal{S}^c$ such that $N^P_{\mathrm{hypo}(V)}(x,V(x))$ is pointed, it holds
\bel{H-p=0}
H(x,\xi)~=~0\qquad\forall \xi\in \partial^{P}V(x),
\eeq
\bel{H-p=00}
\sup_{\omega\in U}\big\{-\xi\cdot f(x,\omega)\big\}~=~0\qquad\forall \xi\in\partial^{\infty}V(x).
\eeq
\end{lemma}
{\bf Proof.} {\bf 1.} Let us first show that 
\[
H(x,\xi)~\geq~0\qquad\forall \xi\in \partial^{P}V(x).
\]
Calling $y(\cdot)$ the optimal trajectory crossing $x$ at time $t>0$,  we have
\[
\begin{split}
-\xi\cdot (y(t-\ve)-x)+ V(y(t-\ve))-V(x)~\leq~\O(1)\cdot \left(|y(t-\ve)-x|^2+|V(y(t-\ve))-V(x)|^2\right)
\end{split}
\]
for all $\ve\in\,]0,t[$. By DPP, we estimate  
\[
{1\over \ve}\cdot \int_{t-\ve}^{t}\Big(\xi\cdot f(y(s), u(s))+r(y(s),u(s))\Big)ds~\leq~\O(1)\cdot\ve,
\]
where $u(\cdot)$ is the control associated to $y(\cdot)$. In particular, this yields
\[
H(x,\xi)~\geq~\limsup_{\ve\to 0+}~{1\over \ve}\cdot \int_{t-\ve}^{t}\Big(-\xi\cdot f(x, u(s))-r(x,u(s))\Big)ds~\geq~\limsup_{\ve\to 0+} \O(1)\cdot \ve~=~0.
\]
Since $N^{P}_{\mathrm{hypo}(V)}(x,V(x))$ is pointed,  by Proposition \ref{Pres-f} the vector $\xi$ can be presented by 
\[
\xi~=~\sum_{i=1}^{d+1}\alpha_i\xi_i+\sum_{i=1}^{d+1}\beta_i \zeta_i
\]
with $\alpha_i,\beta_i\geq 0$ such that $\ds\sum_{i=1}^{d+1}\alpha_i=1$ and 
\[
\xi_i~\in~\partial^{*}V(x) ,\qquad \zeta_i\in\partial^{*,\infty}V(x)\qquad\forall i\in \{1,\dots, d+1\}.
\]
Noticing that 
\bel{Nota1}
H(x,\xi_i)~=~0,\qquad \sup_{\omega\in U}~\big\{-\zeta_i\cdot f(x,\omega)\big\}~=~0,
\eeq
we then achieve 
\[
H(x,\xi)~\leq~\sum_{i=1}^{d+1}\big[\alpha_i\cdot H(x,\xi_i)+\beta_i\cdot \sup_{\omega\in U}~\big\{-\zeta_i\cdot f(x,\omega)\big\}\big]~=~0,
\]
and this yields (\ref{H-p=0}).
\medskip

{\bf 2.} Similarly, if  $\xi$ is in $\partial^{\infty}V(x)$ then by DPP
\[
\begin{split}
{1\over \ve}\cdot  \int_{t-\ve}^{t}\xi\cdot f(y(s), u(s))ds&~\leq~\O(1)\cdot{1\over \ve}\cdot \left(|y(t-\ve)-x|^2+|V(y(t-\ve))-V(x)|^2\right)\\
&~\leq~\O(1)\cdot\ve\qquad\forall \ve>0
\end{split}
\]
and this implies  
\[
\sup_{\omega\in U}\big\{-\xi\cdot f(x,\omega)\big\}~\geq~\limsup_{\ve\to 0+}~\left({1\over \ve}\cdot \int_{t-\ve}^{t}-\xi\cdot f(x, u(s))ds\right)~=~0.
\]
Finally, by Proposition \ref{Pres-f}, we can write $\xi=\ds\sum_{i=1}^{d+1}\beta_i\cdot \zeta_i$ with $\beta_i\geq 0$ and 
\[
 \zeta_i\in\partial^{*,\infty}V(x)\qquad\forall i\in \{1,\dots, d+1\}.
\]
Thus, from (\ref{Nota1}), we get 
\[
\sup_{\omega\in U}\big\{-\xi\cdot f(x,\omega)\big\}~\leq~\sum_{i=1}^{d+1}\beta_i\cdot \sup_{\omega\in U}~\big\{-\zeta_i\cdot f(x,\omega)\big\}~=~0,
\]
and the proof is complete.
\endproof

\begin{theorem} Under the hypotheses of Theorem \ref{Main1}, it holds 
\begin{itemize}
\item [(i).] For every $x\in \Sigma_{V}$ with $\mathrm{dim}\big[\partial^{P}V(x)\cup\partial^{\infty}V(x)\big]< d$, there exist a constant $\sigma>0$ and a Lipschitz arc $\bar{x}(\cdot)$ with $\bar{x}(0)=x$ and $\bar{x}(s)\in\Sigma_{V} \setminus\{x\}$ for all $s\in\, ]0,\sigma]$.
\item [(ii).] Assume that  $H$ satisfies {\bf (H6)}. Then, for every $x\in \mathcal{S}^c\backslash \Sigma_{V,\infty}$, the function $V$ is differentiable along every optimal trajectory starting from $x$ except for the initial and terminal points.
\end{itemize}

\end{theorem}
{\bf Proof.} {\bf 1.} By Theorem  \ref{Main1} and Proposition \ref{Pres-f}, we have that the set $\big(\partial^{\infty}V(x)\backslash\{0\}\big)\cup \partial^{P}V(x)$ is non-empty for all $x\in \mathcal{S}^c$. Assume that $x\in \Sigma_{V}$, we consider two cases:
\begin{itemize}

\item  If  $\partial^{\infty}V(x)\backslash\{0\}$ is empty then by Proposition \ref{hypo-upp}, $V$ is semiconcave in $B(x,\sigma_0)$ for some $\sigma_0>0$. Thanks to the same arguments of \cite[Theorem 3.7]{PCC}, we achieve (i).

\item Otherwise, if there exists a unit vector $\xi\in \partial^{\infty}V(x)$ then by using the same argument in the proof of Lemma \ref{Non-lips}, one can find a sequence $x_n$ converging to $x$ such that $V$ is differentiable at $x_n$ and $\ds\lim_{n\to\infty} DV(x_n)=+\infty$.
Following the same argument in the second step of the proof of Theorem \ref{Re-pre}, we can find an optimal trajectory $\bar{x}(\cdot)=y^{x,\bar{u}}(\cdot)$ for $x$ such that there exists $0\neq \xi^*\in N_0{(x^*)}$ with $x^*=\bar{x}(\tau^{x,\bar{u}})$. Let $\bar{p}(\cdot)$ be the solution to 
\[
\dot{p}(t)~=~-D^Tf(\bar{x}(t),\bar{u}(t)) \cdot p(t),\qquad p(\tau^{x,\bar{u}})~=~\xi^*.
\]
By Proposition \ref{Pro-0}, we have  
\[
\bar p(t)~\in~\partial^{\infty}V(\bar{x}(t))\qquad\forall t\in [0,\tau^{x,\bar{u}}],
\]
and this yields $\bar{x}(t)\in \Sigma_V$ for all $t\in [0,\tau^{x,\bar{u}}]$.
\end{itemize}

{\bf 2.}  For every $x\in \mathcal{S}^c\backslash \Sigma_{V,\infty}$, let $y(\cdot)$ be an optimal trajectory steering  $x$ to the target $\mathcal{S}$ in time $\tau$. For every $x\in \mathcal{S}^c\backslash \Sigma_{V,\infty}$, let $y(\cdot)$ be an optimal trajectory steering  $x$ to the target $\mathcal{S}$ in time $\tau$. 
Assume  that there exists $\xi\in \partial^{\infty} V(y(\tau_0))$ with $|\xi|=1$ for some $\tau_0\in ]0,\tau[$. By Proposition \ref{hypo-upp} and Lemma  \ref{Non-lips}, $y(\tau_0)$ is a  non-Lipschitz point of $V$ and there exists an optimal trajectory $y_{\tau_0}(\cdot) $ steering the point $y(\tau_0)$ to $x^*\in\partial\mathcal{S}$ with $N_0(x^*)\neq \varnothing$. Since the concatenation of $y(\cdot)$ and $y_{\tau_0}(\cdot)$ is an optimal trajectory for $x$, one derives from   Proposition \ref{Pro-0} that $\partial^{\infty} V(x)$ contains a non-zero vector and this yields a contradiction. Thus, one has 
\bel{infty-t}
\partial^{\infty} V(y(t))=\{0\}\qquad\forall t\in \,]0,\tau[.
\eeq
In particular, $N^P_{\mathrm{hypo}(V)}(y(t),V(y(t)))$ is pointed. By Lemma \ref{H-optimal}, we have 
\[
H(y(t),\xi)~=~0\qquad\forall \xi\in \partial^{P}V(y(t))
\]
and {\bf (H6)} implies that $\partial^{P}V(y(t))$ is a singleton. Thus, $V$ is differentiable at $y(t)$ for all $t\in \,]0,\tau[$.
\endproof
\v

{\bf Acknowledgements.} This research by Khai T. Nguyen was partially supported by National Science Foundation grant DMS-2154201. He also  would like to warmly thank Prof. Francis Clarke for asking a question related to  the problem  during his dissertation in 2011.

\end{document}